\documentclass[reqno,12pt]{amsart}
\usepackage{amsfonts,amsmath,amssymb,fullpage,graphicx}
\usepackage{verbatim}

\def\frak{\mathfrak}
\theoremstyle{plain} \newtheorem{Thm}{Theorem}[section]
\theoremstyle{plain} \newtheorem{Cor}[Thm]{Corollary}
\theoremstyle{plain} \newtheorem{Prop}[Thm]{Proposition}
\theoremstyle{plain} \newtheorem{Lemma}[Thm]{Lemma}
\theoremstyle{definition} 
\theoremstyle{definition} \newtheorem{Rem}[Thm]{Remark}
\theoremstyle{definition} \newtheorem{Ex}[Thm]{Example}

\newcommand{\thmlist}{
\renewcommand{\theenumi}{\alph{enumi}}
\renewcommand{\labelenumi}{(\theenumi)}}
\renewcommand{\Re}{\mathop{\rm{Re}}}
\renewcommand{\Im}{\mathop{\rm{Im}}}

\newcommand{\inner}[2]{\langle#1,#2\rangle}

\newcommand{\C}{\ensuremath{\mathbb C}}

\newcommand{\R}{\ensuremath{\mathbb R}}
\newcommand{\Z}{\ensuremath{\mathbb Z}}

\renewcommand{\l}{\lambda}
\renewcommand{\a}{\alpha}
\renewcommand{\b}{\beta}

\newcommand{\fa}{\mathfrak{a}}

\newcommand{\fp}{\mathfrak{p}}
\newcommand{\frakg}{\mathfrak{g}}

\newcommand\SL{\mathop{\rm{SL}}}
\newcommand\GL{\mathop{\rm{GL}}}
\newcommand\SU{\mathop{\rm{SU}}}

\newcommand{\wrhob}{\widetilde{\rho}_\b}

\newcommand{\Res}[1]{\ensuremath{\underset{#1}{{\rm Res}\,}}}

\newcommand{\wt}{\widetilde}

\begin{document}

\makeatletter
\title[Ramanujan's Master Theorem]{Ramanujan's Master Theorem\\ for Riemannian symmetric spaces}
\author{Gestur {\'O}lafsson}
\address{Department of Mathematics, Louisiana State University, Baton Rouge, LA 70803}
\email{olafsson@math.lsu.edu}
\author{Angela Pasquale}
\address{Laboratoire de Math\'ematiques et Applications de Metz (UMR CNRS 7122),
Universit\'e de Lorraine --  Metz, F-57045 Metz, France.}
\email{angela.pasquale@univ-lorraine.fr}
\thanks{The research by {\'O}lafsson was supported by DMS-0801010, DMS-1101337, and the University Paul Verlaine of Metz}

\date{}
\subjclass[2000]{Primary 43A85; Secondary 22E45, 33E20}
\keywords{Ramanujan's Master theorem, integral interpolation formulas, symmetric spaces, spherical transform,
spherical Fourier series, holomorphic extension}

\begin{abstract}
Ramanujan's Master theorem states that, under suitable conditions, the Mellin transform of a power series provides an interpolation formula for the
coefficients of this series. Based on the duality of compact and non-compact reductive Riemannian symmetric spaces  inside a common complexification, we prove an analogue of Ramanujan's Master Theorem for the spherical Fourier
transform of a spherical Fourier series. This extends the results proven by Bertram for Riemannian symmetric spaces of rank-one.
\end{abstract}

\maketitle

\section*{Introduction}

\noindent
Ramanujan's Master Theorem asserts that if a function $f(x)$ can be expanded
around $x=0$ in a power series of the form
$$f(x)=\sum_{k=0}^{\infty} (-1)^k a(k) x^k$$
then
\begin{equation}
\label{eq:Ramanujan}
\int_0^{+\infty} x^{-\l-1} f(x)\; dx= -\frac{\pi}{\sin(\pi \l)} \, a(\l)\,.
\end{equation}
Of course, one needs additional assumptions for this formula to hold, as one can directly see from the
example $a(\l)=\sin(\pi \l)$.

Formula (\ref{eq:Ramanujan}) is contained in Ramanujan's First Quarterly Report to the Board of Studies of the University of Madras in 1913. These reports have never been published; see \cite{Berndt}, p. 295.
However, Hardy presents them in his book on Ramanujan's work \cite{Hardy} and provides a rigorous proof of formula (\ref{eq:Ramanujan}) for a natural class of functions $a$ and a natural set of parameters $\l$, by means of the Residue Theorem.

Let $A$, $P$, $\delta$ be real constants so that $A <\pi$ and $0<\delta\leq 1$. Let
$\mathcal H(\delta)=\{\l \in \C: \Re \l> -\delta\}$. The Hardy class $\mathcal H(A,P,\delta)$ consists of all functions $a:\mathcal H(\delta) \to \C$ that are holomorphic on $\mathcal H(\delta)$
and satisfy the growth condition
$$|a(\l)| \leq C e^{-P(\Re \l)+A|\Im \l|}$$
for all $\l \in \mathcal H(\delta)$. Hardy's version of Ramanujan's Master theorem is the following, see \cite{Hardy}, p. 189.

\begin{Thm}[Ramanujan's Master Theorem]
\label{thm:Ramanujan}
Suppose $a \in \mathcal H(A,P,\delta)$. Then:
\begin{enumerate}
\item The power series
\begin{equation} \label{eq:fpowersum}
f(x)=\sum_{k=0}^{\infty} (-1)^k a(k) x^k
\end{equation}
converges for $0<x<e^{P}$ and defines a real analytic function on this domain.
\item Let $0<\sigma<\delta$. For $0<x<e^{P}$ we have
\begin{equation}\label{eq:extensionPhi}
f(x)=\frac{1}{2\pi i} \, \int_{-\sigma-i\infty}^{-\sigma+i\infty} \frac{-\pi}{\sin(\pi \l)} \, a(\l) x^{\l} \; d\l\,.
\end{equation}
The integral on the right hand side of (\ref{eq:extensionPhi}) converges uniformly on compact subsets of $]0,+\infty[$ and is
independent of the choice of $\sigma$.
\item Formula (\ref{eq:Ramanujan}) holds for the extension of $f$ to $]0,+\infty[$ and for all $\l \in \C$ with $0<\Re \l<\delta$.
\end{enumerate}
\end{Thm}

The last part of Theorem \ref{thm:Ramanujan} is obtained from its second part by applying Mellin's inversion formulas
$$\psi(\l)=\int_0^\infty \phi(x) x^{-\l} \, \frac{dx}{x} \,,\qquad \phi(x)=\frac{1}{2\pi i} \, \int_{-\sigma-i\infty}^{-\sigma+i\infty} \psi(\l) x^{\l} \; d\l$$
to
$$\psi(\l)=\frac{-\pi}{\sin(\pi \l)}\,a(\l),  \qquad \phi(x)=f(x)\,.$$

Formula (\ref{eq:Ramanujan}) can be thought of as an interpolation formula allowing us to reconstruct
$a(\l)$ from the discrete set of its values $\{a(k): k \in \Z^+\}$. In particular, the function $a$ must
vanish identically whenever it vanishes on $\Z^+$. As already noticed by Hardy (\cite{Hardy}, p. 188), this fact and the bound $\pi$ for the exponent $A$ are related to the well-known theorem of Carlson.

An equivalent formulation of Ramanujan's Master Theorem is obtained by replacing $a(\l)\in\mathcal H(A,P,\delta)$ with $A(\l)=a(\l)\Gamma(\l+1)$ where $\Gamma$ is the gamma function. Since
$\Gamma(-\l)\Gamma(\l+1)=-\pi/ \sin(\pi\l)$, this converts the power series (\ref{eq:fpowersum}) into
$$f(x)=\sum_{k=0}^\infty (-1)^k \frac{A(k)}{k!}\, x^k$$
and formula (\ref{eq:Ramanujan}) becomes
\begin{equation}
\label{eq:RamanujanGamma}
\int_0^{+\infty} f(x)x^{-\lambda}\, \frac{dx}{x} =\Gamma(-\l)A(\l)\,.
\end{equation}
It holds for the extension of $f$ to $]0,+\infty[$ and for all $\l\in\C$ with $0<\Re\l<\delta$.
In \cite{Hardy}, \S 11.3, formula (\ref{eq:Ramanujan}) is taken as the standard interpolation formula, though (\ref{eq:RamanujanGamma}) often turns out to be more convenient in the applications.

In \cite{Bertram}, Bertram provided the following group theoretic interpretation of (\ref{eq:Ramanujan}). The functions $x^\l$ ($\l \in \C$) are the spherical functions on $X_G=\R^+$ and the $x^k$ ($k \in \Z$) are the spherical functions on the torus $X_U=U(1)$. Both $X_G$ and $X_U$ can be realized as real forms of their complexification $X_\C=\C^*$. Then (\ref{eq:fpowersum}) and (\ref{eq:extensionPhi}) can be seen respectively as the restriction to $X_U$ and $X_G$ of a ``good" holomorphic function $f$ on $X_\C$.
Let $\widetilde{f}$ and $\widehat{f}$ respectively denote the spherical Fourier transforms of $f$ on $X_G$ and
$X_U$. Then
Ramanujan's formula (\ref{eq:Ramanujan}) becomes
$$\widetilde{f}(\l)=\frac{-\pi}{\sin(\pi \l)} \, a(\l)\,, \qquad \widehat{f}(k)=(-1)^k a(k)\,.$$

By replacing the duality between $U(1)$ and $\R^+$ inside $\C^*$ with the duality
between symmetric spaces of the compact type $X_U=U/K$ and of noncompact type $X_G=G/K$ inside their complexification $X_\C=G_\C/K_\C$, Bertram proved an analogue of Ramanujan's Master theorem for semisimple Riemannian symmetric spaces of rank one.

In this paper we extend Bertram's results first to semisimple Riemannian symmetric spaces of arbitrary rank and then to reductive symmetric spaces.  Our main result is Theorem \ref{thm:RamanujanSymm} which is
Ramanujan's Master Theorem for the spherical Fourier transform on semisimple
Riemannian symmetric spaces. The generalization to reductive symmetric spaces  is then obtained by combining the semisimple case and a multivariable extension of the classical (abelian) Theorem \ref{thm:Ramanujan}; see Theorem \ref{thm:RamanujanSymmRed}.

Let $X_U$ and $X_G$ be dual Riemannian symmetric spaces, respectively of the compact and of the noncompact type, and let $X_\C$ be their complexifiction.  
Our theorem deals with spherical Fourier series on $X_U$ of the form
\begin{equation} \label{eq:fFourier-intro}
f(x)=\sum_{\mu\in \Lambda^+} (-1)^{|\mu|} d(\mu) a(\mu+\rho) \psi_\mu(x)\,, 
\end{equation}
in which the coefficients $a(\mu+\rho)$ are obtained from a holomorphic function $a$ belonging to a certain Hardy class $\mathcal H(A,P,\delta)$ associated with the pair $(X_U,X_G)$ and depending on three real parameters $A,P,\delta$. The function $f$ defines a $K$-invariant holomorphic function on a neighborhood of $X_U$ in $X_\C$. It is then shown to extend holomorphically to a neighborhood of $X_G$ in $X_\C$ by means of the inverse spherical Fourier transform:
\begin{equation}\label{eq:extensionf-symm-intro}
f(x)=\int_{\sigma+i\mathfrak a^*} a(\l)b(\l) \varphi_\l(x) \; \frac{d\l}{c(\l)c(-\l)} \,.
\end{equation}
Formulas (\ref{eq:fFourier-intro}) and (\ref{eq:extensionf-symm-intro}) are the symmetric space analogue of Ramanujan's formulas (\ref{eq:fpowersum}) and (\ref{eq:extensionPhi}). 
An interpolation formula extending (\ref{eq:Ramanujan}) is then obtained from  (\ref{eq:extensionf-symm-intro}) using inversion theorems for the spherical Fourier transform.
The function $b$ occurring in (\ref{eq:extensionf-symm-intro}) is a normalizing factor depending only on the pair 
$(X_U,X_G)$ (and not on the Hardy class). It plays the r\^ole of the function $-(2 i \sin(\pi x))^{-1}$ appearing in Ramanujan's formula  (\ref{eq:extensionPhi}). 
We refer the reader to Theorems \ref{thm:RamanujanSymm} and \ref{thm:RamanujanSymmRed}
for more precise statements of our results and for the unexplained notation in (\ref{eq:fFourier-intro}) and (\ref{eq:extensionf-symm-intro}).

As in the classical proof by Hardy, our principal tool to prove Ramanujan's Master Theorem in the semisimple case is the Residue Theorem.  Our methods are a multivariable generalization of those of Bertram, and in fact we follow several hints to the general rank case that one can find in Section 0.6 of \cite{Bertram}. To deal with the task of taking
residues in our multivariable setting, the main idea is to select the fundamental highest restricted weights as a basis of the
space $\frak a^*_\C$ of spectral parameters and then to work in the corresponding coordinates.

Besides \cite{Bertram}, some other articles have dealt with the extension of Ramanujan's Master Theorem to special classes of semisimple or reductive symmetric spaces. Bertram \cite{BertramCR} and Ding, Gross and
Richards \cite{DGR} have proven this theorem for the (reductive non semisimple) case of symmetric cones.
The version considered in \cite{DGR} corresponds to Ramanujan's formula (B) in Hardy's book \cite{Hardy}, p. 186.
Moreover, Ding \cite{Ding} proved a slightly different version of Ramanujan's Master theorem for Hermitian symmetric spaces by a reduction to the case of symmetric cones from \cite{DGR}. Because of this indirect proving method, the resulting theorem has a more complicate structure than the special case of our Theorem \ref{thm:RamanujanSymm}  for the Hermitian symmetric situation.

This article is organized as follows. In Section \ref{section:prelim} we recall the duality of Riemannian symmetric spaces of compact and noncompact type, we introduce the spherical transform in these contexts, and recall Lassalle's results on the holomorphic extension of spherical Fourier series.
In section \ref{section:Ramanujan} we state Ramanujan's Master theorem for Riemannian symmetric spaces in duality.
In the classical version of the Master Theorem, the interpolated coefficients of the power series are multiplied by a normalizing factor $b(\l)$ which is a constant multiple of the function $\sin(\pi \l)^{-1}$. Similarly, a normalizing function $b(\l)$ appears in the statement of Master's Theorem for symmetric spaces.
In Section \ref{section:b} we explain our choice of the function $b$. Its explicit expression is determined in Section \ref{section:PlancherelDensities}, using the relation between the Plancherel densities for Riemannian symmetric spaces in duality. Examples for Riemannian symmetric spaces of rank-one
or with even root multiplicites are presented. Section \ref{section:estimates} collects several estimates which will be needed in the proof of
Ramanujan's Master Theorem. This theorem is proven in Section \ref{section:proofRamanujanSymm}. In Section
\ref{section:RamanujanRed} we extend our theorem to reductive  Riemannian symmetric spaces.
In the final Section \ref{section:further}, we discuss some further possible extensions of Ramanujan's Theorem for symmetric spaces to related settings and present some open problems.

\subsection*{Notation} We shall use the standard notation $\Z$, $\Z^+$, $\R$, $\R^+$ and $\C$ respectively for the integers, the nonnegative integers, the reals, the nonnegative reals, and the complex numbers.

\section{Preliminaries}
\label{section:prelim}

\subsection{Symmetric spaces and their structure}
\label{subsection:symmetricspaces}
Let $X_U=U/K$ be a Riemannian symmetric space of the compact type. Hence $U$ is a compact connected Lie group, 
and there is an involutive automorphism $\tau$ of $U$ so that $U^\tau_0 \subset K \subset U^\tau$. Here $U^\tau=\{u\in U:\tau(u)=u\}$ and
$U^\tau_0$ is its connected component containing the unit element $e$ of $U$. We will assume that $K$ is connected and $U$ is semisimple. (This assumption on $U$ will be removed in Section \ref{section:RamanujanRed}.)  As $U$ is compact, it admits a finite dimensional faithful representation. So we can assume that $U \subset {\rm U}(m)\subset
{\rm \GL}(m,\C)$ for some $m$. Let $\mathfrak u \subset {\rm M}(m,\C)$ denote the Lie algebra of $U$. The automorphism of $\mathfrak u$ which is the differential of $\tau$ will be indicated by the same symbol. Then $\mathfrak u=\mathfrak k \oplus i\mathfrak p$ where
$\mathfrak k=\{X \in \mathfrak u:\tau(X)=X\}$ is the Lie algebra of $K$ and $i\mathfrak p=\{X \in \mathfrak u:\tau(X)=-X\}$.
Let $U_\C$ denote the analytic subgroup of ${\rm \GL}(m,\C)$ of Lie algebra $\mathfrak u_\C=\mathfrak u \oplus i \mathfrak u$.
Set $\mathfrak g=\mathfrak k \oplus \mathfrak p$ and let $G$ be analytic subgroup of $U_\C$ of Lie algebra $\mathfrak g$.
Then $G$ is a connected noncompact semisimple Lie group containing $K$, and $X_G=G/K$ is a symmetric space
of the noncompact type. Since $\mathfrak g_\C=\mathfrak u_\C$, the group $U_\C$ is a complexification of $G$. We will write
$G_\C=U_\C$.
Let $\mathfrak k_\C=\mathfrak k \oplus i \mathfrak k$ and let $K_\C$ be the connected subgroup of $G_\C$ with Lie algebra $\mathfrak k_\C$. Then $K_\C$ is a closed subgroup of $G_\C$. The symmetric spaces
$X_U=U/K$ and $X_G=G/K$ embed in the complex homogeneous space $X_\C=G_\C/K_\C$ as totally real submanifolds.

A maximal abelian subspace $\mathfrak a$ of $\fp$ is called a Cartan subspace.
The dimension of any Cartan subspace is a constant,
called the \emph{real rank} of $G$ and the \emph{rank} of $X_U$ and $X_G$.  Let $\fa^*$
be the (real) dual space of $\mathfrak a$ and let $\fa_\C^*$ be its
complexification.

Let $\Sigma$ be the set of (restricted) roots of the pair $(\frakg,\fa)$.
It consists of all nonzero $\alpha \in \fa^*$ for
which the vector space $\frakg_{\a} := \{X \in \frakg:
\text{$[H,X]=\a(H)X$  for every $H \in \fa$}\}$ is nonzero. The dimension  of $\frakg_{\a}$ is called the
\emph{multiplicity} of the root $\alpha$ and is denoted by $m_{\a}$.
We fix a set $\Sigma^+$ of positive restricted roots. Then $\Sigma$ is the
disjoint union of $\Sigma^+$ and $-\Sigma^+$. Moreover,
$\fa^+:=\{H\in \fa: \text{$\alpha(H)>0$ for all $\alpha\in
\Sigma^+$}\}$ is an open polyhedral cone called the positive Weyl
chamber.

A root $\a\in \Sigma$ is said to be \emph{unmultipliable} if $2\a \notin
\Sigma$. We respectively denote by
$\Sigma_*$ and $\Sigma_*^+:=\Sigma^+ \cap \Sigma_*$ the sets of unmultipliable roots and of positive unmultipliable roots in $\Sigma$.
The half-sum of the positive roots counted with
multiplicities is denoted by $\rho$: hence
\begin{equation}\label{eq:rho}
\rho=\frac{1}{2}\, \sum_{\alpha\in\Sigma^+} m_\alpha \alpha=
\frac{1}{2}\, \sum_{\b\in\Sigma_*^+} \big(\frac{m_{\b/2}}{2}+ m_\b\big)\b\,\,.
\end{equation}
Here we adopt the usual convention that the multiplicity $m_{\b/2}$ is zero if $\b/2$ is not a root.
By classification, $m_{\b/2}$ is always even.

The Cartan-Killing form $B$ defines a
Euclidean structure on the Cartan subspace $\mathfrak a$. We set
$\inner{X}{Y}:=B(X,Y)$.
We extend this inner product to $\mathfrak a^*$  by
duality, that is we set $\inner{\l}{\mu} :=\inner{H_\l}{H_\mu}$ if
$H_\gamma$ is the unique element in $\fa$ such that
$\inner{H_\gamma}{H}=\gamma(H)$ for all $H \in \fa$. The
$\C$-bilinear extension of $\inner{\cdot}{\cdot}$ to $\fa_\C^*$ will
be denoted by the same symbol.
We shall employ the notation
\begin{equation}
\label{eq:la}
\l_\a=\dfrac{\inner{\l}{\a}}{\inner{\a}{\a}}
\end{equation}
for $\l\in \mathfrak a_\C^*$ and $\a\in\mathfrak a^*$ with $\a\neq 0$.
Notice that $2\l_\a=\l_{\a/2}$.

The \emph{Weyl group} $W$ of $\Sigma$ is the finite group of
orthogonal transformations of $\fa$ generated by the reflections
$r_\a$ with $\a\in \Sigma$, where
\begin{equation*}
r_\a(H):=H-2\frac{\a(H)}{\inner{\a}{\a}}H_\a\,, \qquad H \in \fa\,.
\end{equation*}
The Weyl group action extends to $\fa^*$ by duality, and to $\fa_\C$ and $\fa_\C^*$ by complex linearity.

We set $\mathfrak n=\bigoplus_{\a\in\Sigma^+} \frakg_\a$. Let $N=\exp \mathfrak n$ and
$A=\exp \mathfrak a$ be the connected subgroups of $G$ of Lie algebra $\mathfrak n$ and $\mathfrak a$, respectively. The map $(k, a, n) \mapsto kan$ is an analytic diffeomorphism of the product
manifold $K \times A \times N$ onto $G$.
The resulting decomposition $G =KAN$ is the Iwasawa decomposition
of $G$. Thus, for $g \in G$
we have $g=k(g) \exp H(g) n(g)$ for uniquely determined $k(g)\in K$, $H(g) \in \frak a$
and $n(g)\in G$. We will also need the polar decomposition $G=KAK$:
every $g\in G$ can be written in the form $g=k_1 a k_2$ with $k_1,k_2 \in K$ and $a \in A$.
The element $a$ is unique up to $W$-invariance. It is therefore uniquely determined in $\overline{A^+}$
where $A^+=\exp \mathfrak a^+$.

\subsection{Normalization of measures}
We adopt the normalization of measures from \cite{He3}, Ch. II, \S 3.1. In particular, the Haar
measures $dk$ and $du$ on the compact groups $K$ and $U$ are normalized to have total mass 1.
The Haar measures $da$ and $d\l$ on $A$ and $\frak a^*$, respectively, are normalized so that
the Euclidean Fourier transform
\begin{equation} \label{eq:FourierA}
(\mathcal F_Af)(\l):=\int_A f(a) e^{-i\l(\log a)} \; da\,, \qquad \l\in\frak a^*\,,
\end{equation}
of a sufficiently regular function $f:A \to \C$ is inverted by
\begin{equation} \label{eq:invFourierA}
f(a)=\int_{\frak a^*} (\mathcal F_Af)(\l) e^{i\l(\log a)} \; d\l\,, \qquad a\in A\,.
\end{equation}
The Haar measures $dg$ and $dn$ of $G$ and $N$, respectively, are normalized so that
$dg=e^{2\rho(\log a)}\!dk\,da\,dn$.
Moreover, if $L$ is a Lie group and $P$ is a closed subgroup of $L$, with left Haar measures
$dl$ and $dp$, respectively, then the $L$-invariant measure $d(lP)$ on the homogeneous space
$L/P$ (when it exists) is normalized so that
\begin{equation}
\int_L f(l)\;dl=\int_{L/P} \left( \int_P f(lp)\; dp\right)\,d(lP)\,.
\end{equation}
This condition normalizes the $G$-invariant measure $dx=d(gK)$ on $X_G=G/K$ and fixes the $U$-invariant
measure $d(uK)$ on $X_U=U/K$ to have total mass 1.

\subsection{Spherical functions on $X_G$ and $X_U$}
\label{subsection:spherical}

Let $(\pi,V)$ be an irreducible unitary representation of $G$ or of $U$, and let
$$V^K=\{v \in V:\text{$\pi(k)v=v$ for all $k\in K$}\}$$
be the subspace of the $K$-fixed vectors of $V$. The representation $(\pi,V)$ is said to be \emph{$K$-spherical} 
if $V^K\neq \{0\}$. In this case $\dim V^K=1$.

The spherical functions on $G$ are the matrix coefficients of the (non-unitary) spherical principal series representations.  The spherical function on $G$ of spectral parameter $\l\in\mathfrak a_\C^*$ is the $K$-biinvariant function $\varphi_\l:G \to \C$ given by
\begin{equation} \label{eq:varphil}
\varphi_\l(g)=\int_K e^{(\l-\rho)\big(H(gk)\big)}\, dk
\end{equation}
where $H:G\to \mathfrak a$ is the Iwasawa projection defined earlier. If $e$ denotes the unit element of $G$, then $\varphi_\l(e)=1$. Moreover, $\varphi_\l(g^{-1})=\varphi_{-\l}(g)$ for all $\l\in \mathfrak a^*_\C$ and
$g\in G$.
The function $\varphi_\l(g)$ is real analytic in $g \in G$ and $W$-invariant and entire in $\l\in\mathfrak a_\C^*$.

According to Helgason's theorem, the highest restricted weights of the finite-dimensional $K$-spherical representations of $U$ are the dominant restricted weights, that is the elements of the set
\begin{equation}
\label{eq:LambdaK}
\Lambda^+=\{\mu\in\mathfrak a^*: \text{$\mu_\a\in\Z^+$ for all $\a\in\Sigma^+$}\}\,.
\end{equation}
See \cite{He2}, Ch. V, Theorem 4.1.

Let $\Pi=\{\a_1,\dots,\a_l\}$ be a basis of $\mathfrak a^*$ consisting of simple roots in $\Sigma^+$.
For $j=1,\dots,l$ set
\begin{equation}
\label{eq:betaj}
\beta_j=\begin{cases}
\a_j &\text{if $2\a_j\notin \Sigma$}\\
2\a_j &\text{if $2\a_j\in \Sigma$}
\end{cases}\,.
\end{equation}
Then $\Pi_*=\{\beta_1,\dots,\beta_l\}$ is a basis of $\mathfrak a^*$ consisting of simple roots in $\Sigma_*^+$. Define $\omega_1,\dots, \omega_l \in \mathfrak a^*$ by the conditions
\begin{equation}
\label{eq:omegaj}
(\omega_j)_{\beta_k}=\dfrac{\inner{\omega_j}{\beta_k}}{\inner{\omega_k}{\beta_k}}=\delta_{jk}\,.
\end{equation}
Let $\mu \in \mathfrak a^*$. Then $\mu \in \Lambda^+$ if any only if
$$\mu=\sum_{j=1}^l \mu_j \omega_j \qquad \text{with $\mu_j \in \Z^+$\,, $j=1,\dots,l$}\,.$$
See \cite{He3}, Ch. 2, Proposition 4.23.

The spherical functions on $U$ are matrix coefficients of the finite-dimensional (non-unitary) spherical representations on $U$.
As $U$ is semisimple, the spherical functions on $U$ are parametrized by a subset
$\Lambda_K^+(U)$ of $\Lambda^+$. If $X_U=U/K$ is simply connected, then $\Lambda_K^+(U)=\Lambda^+$.
See \cite{Take}, Ch. II, Theorem 8.2 and Corollary 1, Theorem 6.1(2), and p. 109 and 104.
If $\mu \in \Lambda_K^+(U)$, the corresponding spherical function is given by
\begin{equation} \label{eq:psimu}
\psi_ \mu(u)=\inner{\pi_\mu(u)e_\mu}{e_\mu}\,,
\end{equation}
where $\inner{\cdot}{\cdot}$ denotes the inner product in the space $V_\mu$ of $\pi_\mu$ for which this
representation is unitary and $e_\mu \in V_\mu^K$ is a unit vector.

The spherical functions on $U$ are linked to the
spherical functions on $G$ by holomorphic continuation. More precisely, for every
$\mu \in \Lambda_K^+(U)$, the spherical function $\psi_\mu$ on $U$ with spectral parameter $\mu$ extends
holomorphically to a $K_\C$-biinvariant function on $G_\C$. Its restriction to $G$ is $K$-biinvariant
and coincides with the spherical function $\varphi_{\mu+\rho}$ of $G$.
See \cite{He2}, pp. 540--541, or \cite{BOP}, Lemma 2.5.
The spherical function $\varphi_\l$ of $G$ extends holomorphically to a $K_\C$-biinvariant function on $G_\C$ if and only if $\l$ belongs to the $W$-orbit of $\Lambda_K^+(U)+\rho$.
Notice that, as a matrix coefficient of a unitary representation, $\overline{\psi_\mu(u)}=\psi_\mu(u^{-1})$ for all $u\in U$.

If $L$ is a group acting on a space $X$ and $F(X)$ is a space of functions on $X$, then we shall denote by $F(X)^L$ the subspace of $L$-invariant elements of $F(X)$. In the following, the $K$-biinvariant functions on $G$ (resp. on $U$) will be often identified with the $K$-invariant functions on $X_G=G/K$ (resp. on $X_U=U/K$). In particular, one can consider the $\varphi_\l$'s as $K$-invariant functions on $X_G$ and the $\psi_\mu$ as $K$-invariant functions on $X_U=U/K$. In this case, one can think of the  the spherical functions on $X_U$ as the restrictions of the holomorphic extension of the corresponding spherical functions on $X_G$:
\begin{equation}
\label{eq:varphi-psi}
\psi_\mu=\varphi_{\mu+\rho}|_{X_U}\,, \qquad \mu \in \Lambda_K^+(U)\,.
\end{equation}

\subsection{Spherical harmonic analysis on $X_G$ and $X_U$}
\label{section:spherical-analysis}

The \emph{spherical Fourier transform} of a (sufficiently regular) $K$-invariant function
$f:X_G\to \C$ is the function $\widetilde{f}=\mathcal F_Gf$ defined by
\begin{equation}
\label{eq:sphericalF}
\widetilde{f}(\l)=\mathcal F_Gf(\lambda)=\int_{X_G} f(x) \varphi_{-\l}(x) \; dx
\end{equation}
for all $\l\in \mathfrak a_\C^*$ for which this integral exists.
The Plancherel theorem states that the spherical Fourier transform $\mathcal F_G$
extends as an isometry of $L^2(X_G)^K$ onto $L^2(i\frak a^*, |W|^{-1}|c(\l)|^{-2}d\l)^W$.
Here $|W|$ denotes the order of the Weyl group $W$.
The function $c(\l)$ occurring in the Plancherel density is Harish-Chandra's $c$-function.
It is the meromorphic function on $\frak a_\C^*$ given explicitly by the Gindikin-Karpelevich product formula. In terms of $\Sigma_*^+$, we have
\begin{equation}
\label{eq:c}
c(\l)=c_0 \prod_{\b\in \Sigma_*^+} c_\b(\l)
\end{equation}
where
\begin{equation}
\label{eq:cbeta}
c_\b(\l)= \frac{2^{-2\l_\b} \; \Gamma(2\l_\b)}
{\Gamma\Big(\l_\b+\frac{m_{\b/2}}{4}+\frac{1}{2}\Big)
\Gamma\Big(\l_\b+\frac{m_{\b/2}}{4}+\frac{m_{\b}}{2}\Big)}
\end{equation}
and the constant $c_0$ is given by the condition $c(\rho)=1$.
Formula (\ref{eq:cbeta}) looks slightly different from the usual formula for the $c$-function as found for instance in \cite{He2}, Ch. IV, Theorem 6.4, where it is written in terms of positive indivisible roots ($\a\in\Sigma^+$ with $\alpha/2 \notin \Sigma^+$) rather than in terms of
positive unmultipliable roots.

The spherical Fourier transform has the following inversion formula, which holds for instance a.e. if $f \in L^p(X_G)^K$, with $1 \leq p <2$, and $\widetilde{f} \in L^1(i\mathfrak a^*, |c(\l)|^{-2}d\l)^W$: for almost all $x \in X_G$ we have
\begin{equation}
\label{eq:inversionsphericalG}
f(x)=\frac{1}{|W|} \,
\int_{i\mathfrak a^*} \widetilde{f}(\l) \varphi_\l(x) \; \frac{d\l}{c(\l)c(-\l)}\,.
\end{equation}
See \cite{StantonTomas}, Theorem 3.3.

We shall also need some properties of the $K$-invariant $L^p$-Schwartz spaces
$\mathcal S^p(X_G)^K$ on $X_G$.
Let $1< p\leq 2$ and let $\mathfrak U(\mathfrak g)$ be the universal enveloping algebra of $\mathfrak g_\C$. 
The \emph{$K$-invariant $L^p$-Schwartz space} $\mathcal S^p(X_G)^K$ is the space of all $C^\infty$
$K$-bi-invariant functions $f:G\to \C$ such that for any $D \in \mathfrak U(\mathfrak g)$ and any integer $N \geq 0$ we have
$$
\sup_{g\in G} \big(1+\sigma(g)\big)^N \varphi_0(g)^{-2/p} |(Df)(g)| <\infty\,.
$$
Here $\sigma(g)=\|H\|$ if $g=k_1 \exp(H) k_2$ for $k_1,k_2 \in K$ and $H \in \mathfrak a^+$.
Then $\mathcal S^p(X_G)^K \subset \mathcal S^2(X_G)^K$. Moreover, by identifying as usual $K$-biinvariant functions on $G$ with $K$-invariant functions on $X_G$, we have $\mathcal S^p(X_G)^K \subset L^p(X_G)^K$.

Let $1 < p < 2$ and let $\varepsilon = \frac{2}{p}-1$.
Let ${\rm C}(\varepsilon \rho)^0$ denote the interior of the convex hull of the $W$-orbit of $\varepsilon \rho$ in $\mathfrak a^*$ and let $T_\varepsilon={\rm C}(\varepsilon \rho)^0+
i\mathfrak a^*$ be the tube domain in $\mathfrak a_\C^*$ of base ${\rm C}(\varepsilon \rho)^0$.
The $W$-invariant Schwartz space $\mathcal S(\mathfrak a_\varepsilon^*)^W$ consists of the $W$-invariant holomorphic functions $F: T_\varepsilon \to \C$ such that for every
$u\in S(\mathfrak a_\C)$ and every integer $N\geq 0$ we have
$$\sup_{\l \in T_\varepsilon} (1+\|\l\|)^N |\partial(u)F(\l)| <\infty\,.$$
Then $\mathcal S(\mathfrak a_\varepsilon^*)^W \subset (L^1 \cap L^2)(i\mathfrak a^*,|W|^{-1}
|c(\l)|^{-1} \,d\l)^W$. Moreover, the spherical Fourier transform $\mathcal F_G$ is a bijection of $\mathcal S^p(X_G)^K$ onto the $W$-invariant Schwartz space $\mathcal S(\mathfrak a_\varepsilon^*)^W$.
We refer the reader to \cite{GV}, Ch. 7, and \cite{AnkerSp} for additional information.

The \emph{compact spherical Fourier transform} $\widehat{h}=\mathcal F_U h$ of a (sufficiently regular) $K$-invariant function $h:X_U\to \C$
is usually defined by integration against the spherical functions $\psi_\mu$ on $X_U$.
Because of the relation (\ref{eq:varphi-psi}), we consider $\widehat{h}=\mathcal F_U h$ as the function defined for $\l\in \Lambda_K^+(U) + \rho$ by
\begin{equation}
\label{eq:sphericalU}
\widehat{h}(\l)=\mathcal F_U h(\l)=
\int_{X_U} h(y) \overline{\varphi_{\l}(y)} \; dy
=\int_{X_U} h(y) \varphi_{-\l}(y) \; dy\,.
\end{equation}
The \emph{spherical Fourier series} of $h$ is the formal series on $X_U$ given by
\begin{equation}
\label{eq:sphericalseries}
\sum_{\mu \in \Lambda_K^+(U)} d(\mu) \widehat{h}(\mu+\rho) \psi_{\mu}
=\sum_{\mu \in \Lambda_K^+(U)} d(\mu) \widehat{h}(\mu+\rho) \varphi_{\mu+\rho}\,.
\end{equation}
In (\ref{eq:sphericalseries}), $d(\mu)$ denotes the dimension of the finite dimensional spherical
representation $\pi_\mu$ of $U$ of highest weight $\mu$.
According to Weyl dimension's formula, the function $d$ is a polynomial function on $\mathfrak a_\C^*$.
If $h \in L^2(X_U)^K$, then this series converges to $h$ in $L^2$-norm.
The convergence is absolute and uniform, if $h$ is smooth.
The Plancherel theorem states that compact spherical Fourier
transform $\mathcal F_U$ extends to an isometry of $L^2(X_U)^K$ onto $L^2(\Lambda_K^+(U), d(\mu) d\mu)$
where $d\mu$ is the counting measure.

The spherical harmonic analysis on a general semisimple space $U/K$ of the compact type can be reduced to the simply connected case.
In fact, let $G/K$ and $U/K$ be Riemannian symmetric spaces in duality, as in Section \ref{subsection:symmetricspaces}.
Suppose in that $U/K$ is not simply connected. Let $\wt U$ denote the connected simply connected Lie group with Lie algebra $\mathfrak u$. Let $\wt\theta$ be the involution on $\wt{U}$
with differential equal to the differential of the involution on $U$ associated with $U/K$.
The subgroup $\wt{K}$ of fixed points of $\wt\theta$ is connected. Hence $\wt{U}/\wt{K}$ is a simply connected symmetric space of the compact type. Moreover, there is a subgroup $S$ of the center of $\wt{U}$ so that $U=\wt{U}/S$ and $K=K^*/S$ where $K^*$ is a $\wt\theta$-invariant subgroup of
$\wt{U}$ satisfying $\wt{K}S \subset K^* \subset K_S=\{u \in \wt{U}: u^{-1}\wt\theta(u) \in S\}$. The space $U/K=\wt{U}/K^*$ is then covered by $\wt{U}/\wt{K}$. See \cite{He1}, Ch. VII, Theorem 9.1 and Corollary 9.3.

The group $\wt U$ is connected, simply connected, compact and semisimple. Hence its universal complexification is a connected, simply connected, semisimple, complex Lie group $\wt U_\C$
of Lie algebra $\mathfrak u_\C=\mathfrak g_\C$. Here, as in Section \ref{subsection:symmetricspaces}, we we have set $\mathfrak g=\mathfrak k \oplus \mathfrak p$ and $\mathfrak u=\mathfrak k\oplus i\mathfrak p$ for the Lie algebras of $G$ and $U$, respectively.

Let $G^\sharp$ be the connected Lie subgroup of $\wt U_\C$ with Lie algebra $\mathfrak g$.
Then the inclusion of $G^\sharp$ in $\wt U_\C$ gives the universal complexification of
$G^\sharp$. Moreover, $G^\sharp$ has finite center. Since $\mathfrak k \subset \mathfrak g$ and $\wt K$ is connected, we also have $\wt K \subset G^\sharp$. Hence $\wt K$ is a maximal compact subgroup of $G^\sharp$. Thus, $G^\sharp/\wt K$ is a symmetric space of the noncompact type, $G^\sharp$ is a connected subgroup of its universal complexification
$\wt U_\C$, and the compact dual of $G^\sharp/\wt K$ is $\wt U/\wt K$, where $\wt U$ is the connected subgroup of $\wt U_\C$ of Lie algebra $\mathfrak u$. This replaces the original
dual pair $(G/K,U/K)$ with the dual pair $(G^\sharp/\wt K,\wt U/\wt K)$ with $\wt U/\wt K$
simply connected.

Let $\pi: \wt{U}/\wt{K} \to \wt{U}/K^*=U/K$ be the covering map. By composing with $\pi$, we can identify a $K$-invariant function $f$ on $U/K$ with a $K^*$-invariant function on $ \wt{U}/\wt{K}$. The space of $K$-invariant functions on $U/K$ can then be considered as the subspace of $\wt{K}$-invariant functions on $\wt{U}/\wt{K}$ that, moreover, are $K^*$-invariant. The compact spherical Fourier transform of $f$ as a $K$-invariant function on $U/K$ is then identified with the restriction to
$\Lambda_K^+(U)$ of the compact spherical Fourier transform of $f$ as $K^*$-invariant function on $\wt{U}/\wt{K}$. On the noncompact side, the canonical isomorphism of $G/K$ and
$G^\sharp/\wt K$ allows us to identify $K$-invariant functions on $G/K$ with
$\wt K$-invariant functions on $G^\sharp/\wt K$. Under these identifications, the noncompact spherical transforms on these two symmetric spaces agree.

\subsection{Holomorphic extension of spherical Fourier series}
Set $\mathfrak b:=i\mathfrak a$, and let $\|\cdot\|$ be a $W$-invariant norm on $\mathfrak b$. Endow $\mathfrak a^*$ with the dual norm, still denoted by the same symbol. Let $B=\{H \in \mathfrak b: \|H\| < 1\}$ be the open unit ball in $\mathfrak b$. Let $o_\C=eK_\C$ be the base point in $X_\C=G_\C/K_\C$.
Under suitable exponential decay of their coefficients, Lassalle proved in \cite{Lassalle} the normal convergence of the spherical Fourier series on  $U$-invariant domains in $X_\C$ of the form $D_\varepsilon=U \exp(i\varepsilon B)\cdot o_\C$.
Notice that $D_\varepsilon$ is a neighborhood of $X_U$ in $X_\C$ and that
$D_\varepsilon \supset K\exp(i\varepsilon B)\cdot o_\C$. Moreover, $K\exp(i\varepsilon B)\cdot o_\C$ is an open neighborhood of $o_\C$ in $KA \cdot o_\C$, which is the image of $X_G$ in its embedding in $X_\C$.
For the reader's convenience, we collect the results which will be needed in the following.

\begin{Thm}
\label{thm:Lassalle}
\begin{enumerate}
\thmlist
\item Let $F:\Lambda^+_K(U)\to \mathbb{C}$.
Suppose there are constants $C>0$ and $\varepsilon >0$ so that for all $\mu \in \Lambda^+_K(U)$
$$|F(\mu)|\leq C e^{-\varepsilon \|\mu\|}\,.$$
Then the spherical Fourier series
$$\sum_{\mu \in \Lambda^+_K(U)} d(\mu) F(\mu) \psi_\mu(x)$$
converges normally on compact subsets of $D_\varepsilon$. Its sum is therefore a holomorphic $K$-invariant function on $D_\varepsilon$.
\item
Conversely, suppose that $h$ is a continuous $K$-invariant function on $X_U$ admitting a holomorphic extension to a neighborhood of $X_U$ in $X_\C$. Let $F(\mu)=\widehat{h}(\mu+\rho)$ be the Fourier coefficients of $h$. Then there are constants $C>0$ and $\varepsilon >0$ so that for all $\mu \in \Lambda^+_K(U)$

$$|F(\mu)|\leq C e^{-\varepsilon \|\mu\|}\,.$$
\end{enumerate}
\end{Thm}
\begin{proof}
For the proof of part (a), we follow \cite{Lassalle}, p. 189.
It is enough to prove the normal convergence of the series on compact sets of the
form $\overline{D_r}=U\exp(ir\overline{B}).o_\C$ where $0<r<\epsilon$ and
$\overline{B}=\{H\in\mathfrak b:\|H\|\leq 1\}$ is the closed unit ball in $\mathfrak b$. Set $\mathfrak b^+=i\mathfrak a^+$.
If $x \in X_\C$, then we can write $x=u\exp\big(iA(x)\big) \cdot o_\C$ for a unique $A(x) \in \overline{\mathfrak b^+}$ and some (non unique) $u \in U$; see \cite{Lassalle}, Th\'eor\`eme 1, p. 177.
If $A(x) \in r\overline{B}$, then, by \cite{Lassalle}, Proposition 12, p. 184, one has
\begin{equation*}
|\psi_\mu(x)| \leq e^{\mu\big( i A(x)\big)} \leq
e^{\sup_{H\in r\overline{B}} \mu(i H)}= e^{r\|\mu\|}\,.
\end{equation*}
Thus
$$|d(\mu) F(\mu) \psi_\mu(x)|\leq d(\mu) e^{-(\varepsilon -r)\|\mu\|}\,,$$
which implies the convergence of the series as $d(\mu)$ is a polynomial in $\mu$.

Part (b) is a special case of Proposition V.2.3 in \cite{Faraut}. It is proven using Cauchy's inequalities for the Fourier coefficients of $h$.
\end{proof}

We observe that Theorem \ref{thm:Lassalle} holds in the general case where $U$ is reductive. This will be needed in Section
\ref{section:RamanujanRed}.

\subsection{Coordinates in $\mathfrak a_\C^*$ and tubes domains around $i\mathfrak a^*$}
\label{subsection:coordinates}
We choose $\Pi_*=\{\omega_1,\dots,\omega_l\}$ as basis of $\mathfrak a^*$. For $\l\in\mathfrak a^*$ we have
\begin{equation}\label{eq:lambdaSigmaast}
\l=\sum_{j=1}^l \l_j \omega_j \qquad\text{with}\qquad \l_j:=\l_{\b_j}=\frac{\inner{\l}{\b_j}}{\inner{\b_j}{\b_j}} \,.
\end{equation}
Set $\mathfrak a^*_+=\{\l\in \mathfrak a^* : \text{$\l_\b\geq 0$ for all $\b\in\Sigma^+_*$}\}$.
By identifying $\l\equiv (\l_1,\dots,\l_l)$, we obtain the correspondences $\mathfrak a^* \equiv \R^l$,
  $\mathfrak a^*_+ \equiv (\overline{\R^+})^l$ and $\Lambda^+ \equiv (\Z^+)^l$.
For $\mu=\sum_{j=1}^l \mu_j \omega_j\in\Lambda^+$, we define
\begin{equation}
\label{eq:heightmu}
|\mu|=\mu_1+\dots+\mu_l\,.
\end{equation}

Set
\begin{equation} \label{eq:defrhoj}
\rho=\sum_{j=1}^l \rho_j \omega_j\,.
\end{equation}
Since $\b_j$ is a multiple of a simple root, we have
\begin{equation}
\rho_j=\frac{1}{2} \left( m_{\b_j}+ \frac{m_{\b_j/2}}{2}\right)\,.
\end{equation}
For an arbitrary $\b \in \Sigma_*^+$, we set
\begin{equation}
\label{eq:wrhob}
\wrhob=\frac{1}{2} \left( m_{\b}+ \frac{m_{\b/2}}{2}\right)\,.
\end{equation}
Notice that $\wrhob=\rho_j=\rho_{\b_j}$ if $\b=\b_j$, but $\wrhob\neq \rho_\b$ in general.

Let $\delta>0$. We consider the following tube domains in
$\mathfrak a^*_\C$ around the imaginary axis:
\begin{align}
\label{eq:Trho}
T_\delta&=\{\l\in\mathfrak a^*_\C: \text{$|\Re\l_\b|<\delta \wrhob$ for all $\b\in\Sigma_*^+$}\}\,,\\
T'_{\delta}&=\{\l \in \mathfrak a_\C^*: \text{$|\Re\l_j|<\delta \rho_j$ for all $j=1,\dots,l$}\}\,, \\
\label{eq:Trhodue}
T''_{\delta}&=\{\l \in \mathfrak a_\C^*: \text{$\Re\l_j<\delta \rho_j$ for all $j=1,\dots,l$}\}\,.
\end{align}
In the following we shall denote by $B(T)$ the base in $\mathfrak a^*$ of the tube domain $T$ in $\mathfrak a^*_\C$. Hence
$T=B(T)+i\mathfrak a^*$ and $B(T)=T \cap \mathfrak a^*$. The following lemma is rather standard; we provide a proof
for the sake of completeness.

\begin{Lemma}
\label{lemma:Trho}
Let $w_0$ be the longest element of $W$. Then
\begin{equation}\label{eq:Trhow0}
T'_\delta=T''_\delta \cap w_0(T''_\delta)
\end{equation}
and
\begin{equation} \label{eq:Trhouno}
T_{\delta}=\bigcap_{w \in W} w(T'_\delta)=\bigcap_{w \in W} w(T''_\delta)\,.
\end{equation}
In particular, $T_\delta$ is the largest $W$-invariant tube domain contained in $T'_{\delta}$.
Moreover,
\begin{equation} \label{eq:TrhoCrho}
T_\delta={\rm C}(\delta \rho)^0+i\mathfrak a^*
\end{equation}
where ${\rm C}(\nu)$ is the the convex hull of the $W$-orbit $\{w\nu: w \in W\}$ of $\nu\in \mathfrak a^*$
and ${\rm C}(\nu)^0$ is its interior.
\end{Lemma}
\begin{proof}
It is enough to prove that the bases in $\mathfrak a^*$ of the considered tube domains are equal.

The element $w_0$ maps $\Sigma^+$ into $-\Sigma^+$. Hence $w_0\rho=-\rho$. Moreover, $w_0$ maps the basis $\Pi$ of simple roots in $\Sigma^+$ into the basis $-\Pi$ in $-\Sigma^+$. Hence
it maps $\{\b_1,\dots,\b_l\}$ into $\{-\b_1,\dots,-\b_l\}$. Furthermore, it satisfies $w_0^{-1}=w_0$. Therefore
\begin{align*}
w_0\big(B(T''_\delta)\big)
&=\{w_0\l\in \mathfrak a^*:\text{$\inner{\l-\delta \rho}{\b_j}<0$ for all $j=1,\dots,l$} \}\\
&=\{w_0\l\in \mathfrak a^*:\text{$\inner{\l-\delta \rho}{w_0\b_j}>0$ for all $j=1,\dots,l$} \}\\
&=\{w_0\l\in \mathfrak a^*:\text{$\inner{w_0\l-\delta w_0\rho}{\b_j}>0$ for all $j=1,\dots,l$} \}\\
&=\{\l\in \mathfrak a^*:\text{$\inner{\l+\delta \rho}{\b_j}>0$ for all $j=1,\dots,l$} \}\\
&=\{\l\in \mathfrak a^*:\text{$\l_j>-\delta \rho_j$ for all $j=1,\dots,l$} \}\,.
\end{align*}
This proves (\ref{eq:Trhow0}).

Set $C=B\big( \bigcap_{w \in W} w(T'_\delta)\big)=\bigcap_{w \in W} w\big(B(T'_\delta)\big)$.
We have $B(T_\delta) \subset B(T'_\delta)$ as $\widetilde{\rho}_{\b_j}=\rho_j$. Hence $B(T_\delta) \subset C$ as $B(T_\delta)$ is $W$-invariant.
Conversely, let $\l\in C$ and let $\b \in \Sigma_*$. Since $\Sigma_*$ is a reduced root system and the corresponding Weyl group is $W$, there is $j\in\{1,\dots,l\}$ and $w \in W$ so that $w\b=\b_j$.
Since $w\l\in C \subset B(T'_\delta)$, we obtain that $|\l_\b|=|(w\l)_{w\b}|=|(w\l)_{\b_j}|<\delta \rho_j=\delta \wrhob$\,.
Thus $\l \in B(T_\rho)$. The second equality in (\ref{eq:Trhouno}) is a consequence of (\ref{eq:Trhow0}).

For the final result, we have by \cite{He2}, Lemma 8.3 (i), that $C(\delta\rho)=\bigcap_{w \in W} w\big(\delta\rho-\mathfrak a_+^*\big)$. Notice that
$$\delta\rho-\mathfrak a_+^*=\{\delta\rho-\l\in\mathfrak a^*:\text{$\l_j \leq 0$ for all $j=1,\dots,l$}\}
=\{\l\in\mathfrak a^*: \text{$\l_j\leq \delta\rho_j$ for all $j=1,\dots,l$}\}\,.$$
Thus, by (\ref{eq:Trhouno}),
$$C(\delta\rho)^0=\bigcap_{w \in W} w\big((\delta\rho-\mathfrak a_+^*)^0\big)=
\bigcap_{w \in W} w\big(B(T''_\delta)\big)\,=B(T_\delta)\,.$$
\end{proof}
By (\ref{eq:TrhoCrho}) and a theorem by Helgason and Johnson, $T_1=C(\rho)^0+i\mathfrak a^*$ is the interior of the set of parameters $\l\in\mathfrak a_\C^*$ for which the spherical function $\varphi_\l$ is bounded. See \cite{He2}, Theorem 8.1.

\section{Statement of Ramanujan's Master theorem for symmetric spaces}
\label{section:Ramanujan}
\noindent
Let $X_G=G/K$ and $X_U=U/K$ be semisimple symmetric spaces in duality inside their common complexification
$X_\C$ as in Section \ref{subsection:symmetricspaces}. We suppose that $X_U$ is simply connected.

Let $A$, $P$, $\delta$ be constants so that $A <\pi$, $P>0$ and
$0<\delta\leq 1$.
Let
\begin{equation}\label{eq:Hdelta}
\mathcal H(\delta)=\{\l \in \mathfrak a_\C^*: \text{$\Re \l_\b > -\delta\,\wt\rho_\b$ for all $\b\in\Sigma_*^+$}\}\,.
\end{equation}
The \emph{Hardy class} $\mathcal H(A,P,\delta)$ consists of the functions $a:\mathcal H(\delta) \to \C$ that are holomorphic on $\mathcal H(\delta)$
and so that
\begin{equation}
\label{eq:defHardyclass}
|a(\l)| \leq C \prod_{j=1}^l e^{-P(\Re \l_j)+A|\Im \l_j|}
\end{equation}
for some constant $C \geq 0$ and for all $\l \in \mathcal H(\delta)$.

We denote by $\|X\|$ the norm of $X\in \mathfrak a$ with respect to the $W$-invariant norm which is induced by the Killing form.
The same notation is also employed for the corresponding norm on $\mathfrak a^*$.
Recall from Section \ref{section:spherical-analysis} the notation $d$ for the polynomial function on $\mathfrak a_\C^*$
given by  Weyl dimension formula. An explicit expression of $d$ in terms of the positive unmultipliable roots
will be given in formula (\ref{eq:d-poly}).

\begin{Thm}[Ramanujan's Master Theorem for semisimple Riemannian symmetric spaces]
\label{thm:RamanujanSymm}

Let $b$ be the meromorphic function on $\mathfrak a_\C^*$ defined by the equality
\begin{equation}
\label{eq:b}
\frac{b(\l)}{c(\l)c(-\l)}=\left(\frac{i}{2}\right)^l d(\l-\rho)\; \prod_{j=1}^l \frac{1}{\sin\big(\pi(\l_j-\rho_j)\big)}
\end{equation}
and let
\begin{multline}
\label{eq:TSigmam}
T_{\Sigma,m}=\{\l\in\mathfrak a_\C^*: \text{$|\Re\l_\b|<1$ for all $\b\in\Sigma_*^+$ with $(m_{\beta/2})/2$ even}\}  \\ \cap \{\l\in\mathfrak a_\C^*: \text{$|\Re\l_\b|<1/2$ for all $\b\in\Sigma_*^+$ with $(m_{\beta/2})/2$ odd}\}\,.
\end{multline}
Suppose $a \in \mathcal H(A,P,\delta)$.
\begin{enumerate}
\item
Set
\begin{equation}
\label{eq:Omega}
\Omega=\max_{j=1,\dots,l} \|\omega_j\|\,.
\end{equation}
Then the  spherical Fourier series
\begin{equation} \label{eq:fFourier}
f(x)=\sum_{\mu\in \Lambda^+} (-1)^{|\mu|} d(\mu) a(\mu+\rho) \psi_\mu(x)
\end{equation}
converges normally on  compact subsets of $D_{P/\Omega}=U \exp\big((P/\Omega)B\big) \cdot o_\C$
where $B=\{H \in \mathfrak a: \|H\|<1\}$ is the open unit ball in $\mathfrak a$. Its sum is a $K$-invariant holomorphic function on the neighborhood $D_{P/\Omega}$ of $X_U$ in $X_\C$.

\item Let $T_\delta$ be the tube domain in (\ref{eq:Trho}) and let
$\sigma\in B(T_\delta)$. Then for $x=\exp H \in A$ with $\|H\|<P/\Omega$, we have
\begin{equation}\label{eq:extensionf-symm}
f(x)=\frac{1}{|W|} \int_{\sigma+i\mathfrak a^*} \left( \sum_{w \in W} a(w\l)b(w\l)\right) \varphi_\l(x) \; \frac{d\l}{c(\l)c(-\l)} \,.
\end{equation}
The integral on the right hand side of (\ref{eq:extensionf-symm}) is independent of the choice of $\sigma$. It converges uniformly on compact subsets of $A$ and extends to a holomorphic $K$-invariant function on a neighborhood of $X_G$ in $X_\C$.
\item The extension of $f$ to $X_G$
satisfies
$$\int_{X_G} |f(x)|^2 \, dx =\frac{1}{|W|} \int_{i\mathfrak a^*} \Big| \sum_{w \in W} a(w\l)b(w\l)\Big|^2 \; \frac{d\l}{|c(\l)|^2}
\,.$$
Moreover,
\begin{equation} \label{eq:Ramanujansymm}
\int_{X_G} f(x)\varphi_{-\l}(x)\; dx=  \sum_{w \in W} a(w\l)b(w\l)
\end{equation}
for all $\l \in T_\delta \cap T_{\Sigma,m}$. More precisely, the integral on the left-hand side of (\ref{eq:Ramanujansymm}) converges in $L^2$-sense and absolutely on $i\mathfrak a^*$.
It defines a $W$-invariant holomorphic function on a $W$-invariant tube domain around $i\mathfrak a^*$, and (\ref{eq:Ramanujansymm}) extends as an identity between holomorphic functions on $T_\delta \cap T_{\Sigma,m}$.
\end{enumerate}
\end{Thm}

\begin{Rem}
The function $b$ occurring in Theorem \ref{thm:RamanujanSymm} plays the role of the function $\frac{-1}{2\pi i}\,
\frac{\pi}{\sin(\pi x)}=\frac{i}{2}\, \frac{1}{\sin(\pi x)}$ of the classical formula by Ramanujan. A more explicit formula for this function will be given in Corollary \ref{cor:b}. The tube domain $T_{\Sigma,m}$ is linked to the singularities of $b$. Indeed,
$T_{\Sigma,m} \cap T_\delta$ is the largest domain on which the function $\sum_{w\in W} a(w\l)b(w\l)$ is holomorphic for any $a \in \mathcal H(A,P,\delta)$. See Lemma \ref{lemma:Schwartztube}. Furthermore, the condition $0<\delta\leq 1$ guarantees the convergence of the integral
in (\ref{eq:extensionf-symm}), as in the classical case.
\end{Rem}

\begin{Rem} \label{rem:equivalentnorms}
Set $\|\l\|_1=\sum_{j=1}^l |\l_j|$ for $\l=\sum_{j=1}^l \l_j\omega_j \in \mathfrak a^*$.
Then $\|\cdot\|_1$ is a norm on $\mathfrak a^*$.
The growth condition (\ref{eq:defHardyclass}) can be written as $|a(\l)| \leq C e^{-P\|\Re\l\|_1+
A\|\Im \l\|_1}$. Since the spherical transform maps into $W$-invariant functions, it is sometimes more convenient to consider estimates with respect to the $W$-invariant norm $\|\cdot\|$ on $\mathfrak a^*$  which is associated with the Killing form. This can of course be done by norm equivalence.
We shall denote by $c_1$ and $c_2$ two positive constants so that
$c_1 \|\l\| \leq \|\l\|_1 \leq c_2 \|\l\|$ for all $\l\in \mathfrak a^*$.
By (\ref{eq:lambdaSigmaast}), one can choose for instance $c_1=\Omega^{-1}$ with $\Omega$ as in (\ref{eq:Omega}) and $c_2=\sum_{j=1}^l \|\beta_j\|^{-1}$.
\end{Rem}

\begin{Rem} \label{rem:positiveP}
Hardy's version of Ramanujan's Master Theorem holds for functions in $\mathcal H(A,P,\delta)$ when $P$ is an arbitrary real number, whereas Theorem \ref{thm:RamanujanSymm} is stated only for $P>0$.
This assumption cannot be removed. Notice first that if $a \in \mathcal H(A,P,\delta)$, then
for all $\mu \in \Lambda^+$ we have $|a(\mu+\rho)| \leq C e^{-P\sum_{j=1}^l \mu_j} \leq
C e^{-Pc \|\mu\|}$ where $c$ is a positive constant (see Remark \ref{rem:equivalentnorms}).
According to Theorem \ref{thm:Lassalle} (b), $P>0$ is a necessary condition for the elements $(-1)^{|\mu|} a(\mu+\rho)$ to be the Fourier coefficients of a continuous $K$-invariant function on $X_U$ admitting a holomorphic extension to a neighborhood of $X_U$ in $X_\C$. In this case, the Fourier series converges in an open domain of $X_G$ containing the base point $o$. It is natural to ask whether a condition $P\geq 0$ could still allow the absolute convergence of the series $\sum_{\mu\in\Lambda^+}  (-1)^{|\mu|} d(\mu) a(\mu+\rho) \psi_\mu(x)$ for $x$ in some open domain in $X_G$. The answer is negative. Indeed, by $K$-biinvariance, we can restrict ourselves to domains in $\overline{A^+}\equiv \overline{A^+}\cdot o$.
According to Proposition IV.5.2 in \cite{Faraut}, one has the estimate
\begin{equation*}
\psi_\mu(\exp H) \geq c(\mu+\rho) e^{\mu(H)}
\end{equation*}
for $\exp H \in \overline{A^+}$ and $\mu \in \Lambda^+$. Formulas (\ref{eq:c}) and (\ref{eq:cbeta}) show that $c(\mu+\rho)>0$ for $\mu \in \Lambda^+$.
If $\sum_{\mu\in\Lambda^+}  (-1)^{|\mu|} d(\mu) a(\mu+\rho) \psi_\mu(\exp H)$ converges absolutely, then there is a constant $C_H>0$ so that
$$|a(\mu+\rho)| c(\mu+\rho) e^{\mu(H)} \leq |a(\mu+\rho) \psi_\mu(\exp H)| \leq C_H$$
for all $\mu\in\Lambda^+$. Hence $|a(\mu+\rho)|\leq  c(\mu+\rho)^{-1} e^{-\mu(H)}$ has exponential decay in $\mu$ for $H\neq 0$.

Notice also that in the classical version of Ramanujan's Master Theorem, the series (\ref{eq:fpowersum}) converges normally on compact subsets of the domain $|z|<e^P$. This set is a neighborhood of $U(1)=\{z \in \C:|z|=1\}$ if and only if $P>0$.
The restriction $P>0$ is therefore intrinsically related to the geometric interpretation of Ramanujan's Master Theorem as a relation of mathematical objects on symmetric spaces in duality.
\end{Rem}

\begin{Rem}
A slightly more general definition of Hardy class could be obtained by replacing the constants $\delta$ and $P$ respectively by a
$W$-invariant multiplicity function $\delta:\Sigma_*\to ]0,1]$ and a vector $P=(P_1,\dots, P_l)$ with $P_j>0$ for all $j$.
The version of Ramanujan's Master Theorem for this notion of Hardy class would not present any new difficulty with respect to the
version proven here. We have nevertherless preferred to present the case where $\delta$ and $P$ are constants, to keep the
notation as simple as possible.
\end{Rem}

\begin{Rem}
\label{rem:notsimplyconnected}
Suppose that $X_G=G/K$ and $X_U=U/K$ is a pair of symmetric spaces in duality, as in Section \ref{subsection:symmetricspaces}, with
$U/K$ not simply connected. As in Section \ref{subsection:spherical}, we can identify the $K$-invariant functions $f$ on $U/K$ with the $K^*$-invariant functions on its simply connected cover $\wt{U}/\wt{K}$. In this way, the compact spherical Fourier transform of $f$ is the restriction to $\Lambda_K^+(U) \subset \Lambda^+$ of the transform of $f$ as a function on $\wt{U}/\wt{K}$. Likewise, we identify the $K$-invariant functions on $G/K$ with the $\wt K$-invariant functions on $G^\sharp/\wt K$. See Section \ref{subsection:spherical}.
Ramanujan's Master theorem for $X_G$ and $X_U$ is then obtained from the one for $\wt{U}/\wt{K}$ and $G^\sharp/\wt K$ by replacing the Hardy class $\mathcal H(A,P,\delta)$ with its subspace
$$\mathcal H_{U/K}(A,P,\delta)=\{a\in \mathcal H(A,P,\delta): \text{$a(\mu)=0$ for all $\mu \in \Lambda^+\setminus \Lambda^+_K(U)$}\}\,.$$
\end{Rem}
\bigskip

As in the classical case, there is an equivalent formulation of Ramanujan's Master Theorem using the gamma function.
For simplicity, we only consider the case when $X_U$ is simply connected. The general case can be dealt with as in Remark \ref{rem:notsimplyconnected}. Let $B(\l)$ be the meromorphic function on $\mathfrak a_\C^*$ defined by the equality
$$\frac{B(\l)}{c(\l)c(-\l)}=\left(\frac{1}{2\pi i}\right)^l d(\l-\rho) \prod_{j=1}^l \Gamma(-\l_j+\rho_j)\,,$$
i.e.
$$B(\l)=b(\l) \prod_{j=1}^l \frac{1}{\Gamma(\l_j-\rho_j+1)}\,.$$
Replace $a(\l) \in \mathcal H(A,P,\delta)$ with $A(\l)=a(\l) \prod_{j=1}^l \Gamma(\l_j-\rho_j+1)$
in Theorem \ref{thm:RamanujanSymm}. Notice that $A$ is in general no longer holomorphic on $\mathcal H(\delta)$.
However, $A(\lambda)$ is holomorphic on
$\{\l \in \mathfrak a_\C^*:\text{$\l_j> \rho_j-1$ for all $j=1,\dots,l$}\}$ and $A(\l)B(\l)=a(\l)b(\l)$.
The power series (\ref{eq:fFourier}) now becomes
\begin{equation}
\label{eq:FFourier}
F(x)=\sum_{\mu\in\Lambda^+} (-1)^{|\mu|} d(\mu) \frac{A(\mu+\rho)}{\mu!} \, \psi_\mu(x)
\end{equation}
where
\begin{equation}
\label{eq:factorial}
\mu ! =\mu_1 ! \cdots \mu_l !  \qquad\text{for \qquad $\mu=\sum_{j=1}^l \mu_j \omega_j$}\,.
\end{equation}
Moreover, the extension of $F$ to $X_G$ satisfies
\begin{equation}
\label{eq:RamanujanSymmF}
\int_{X_G} F(x)\varphi_{-\l}(x)\; dx=\sum_{w \in W} A(w\l)B(w\l)\,.
\end{equation}

Other formulas can be deduced from (\ref{eq:Ramanujansymm}) or (\ref{eq:RamanujanSymmF}) by formal manipulations. We mention for instance the following, which is the analogue Ramanujan's formula in \cite{Berndt}, Corollary (i), p. 318.

Let $\wt B(\l)$ be the meromorphic function defined on $\mathfrak a_\C^*$ by the equality
$$\frac{\wt B(\l)}{c(\l)c(-\l)}=\left(\frac{1}{4\pi i}\right)^l \pi^{l/2} d(\l-\rho) \prod_{j=1}^l
\frac{\Gamma(-\l_j+\rho_j)}{\Gamma((\l_j-\rho+1)/2)}\,,$$
i.e.
$$\wt B(\l)=b(\l) \prod_{j=1}^l \frac{\cos(\pi(\l_j-\rho_j)/2)}{\Gamma(\l_j-\rho_j+1)}\,.$$
Let
$$
\wt A(\l)= a(\l)\; \frac{\Gamma(\l_j-\rho_j+1)}{\cos(\pi(\l_j-\rho_j)/2)}
$$
with $a(\l) \in \mathcal H(A,P,\delta)$.
The power series (\ref{eq:fFourier}) becomes
\begin{equation}
\label{eq:FtildeFourier}
\wt F(x)=\sum_{\mu\in\Lambda^+} d(2\mu) \frac{\wt A(2\mu+\rho)}{(2\mu)!} \, \psi_{2\mu}(x)
\end{equation}
and the extension of $\wt F$ to $X_G$ satisfies
\begin{equation}
\label{eq:RamanujanSymmFtilde}
\int_{X_G} \wt F(x)\varphi_{-\l}(x)\; dx=\sum_{w \in W} \wt A(w\l)\wt B(w\l)\,.
\end{equation}

\section{The function $b$}
\label{section:b}
\noindent
In this section we explain our choice of the function $b$ occurring in the statement of Ramanujan's Master theorem, Theorem \ref{thm:RamanujanSymm}.

We are looking for a meromorphic function $b:\mathfrak a_\C^* \to \C$ with the following property: For every function $a:\mathfrak a_\C^* \to \C$ in the Hardy class, we have
\begin{equation}
\label{eq:keyequality-b-uno}
\frac{1}{|W|} \int_{i\mathfrak a^*} \Big( \sum_{w \in W} a(w\l)b(w\l)\Big) \varphi_\l(x) \; \frac{d\l}{c(\l)c(-\l)} = \sum_{\mu \in \Lambda^+} (-1)^{|\mu|} d(\mu) a(\mu+\rho) \varphi_{\mu+\rho}(x)
\end{equation}
for all $x \in X_G$ sufficiently close to the base point $o=eK$.

The right-hand side of (\ref{eq:keyequality-b-uno}) will represent the holomorphic extension to a neighborhood of $X_U$ in $X_\C$ of the $K$-invariant function
$f_1:X_U\to \C$ defined by the the spherical Fourier series
$$f_1(y)=\sum_{\mu \in \Lambda^+} (-1)^{|\mu|} d(\mu) a(\mu+\rho) \psi_{\mu}(y)\,, \qquad y \in X_U\,.$$
The right-hand side will give a $K$-invariant function
$f_2:X_G \to \C$ having the $W$-invariant function $\sum_{w \in W} a(w\l)b(w\l)$, $\l \in i\mathfrak a^*$, as noncompact spherical transform.

Observe that, since the spherical function $\varphi_\l$ and the Plancherel density $[c(\l)c(-\l)]^{-1}$ are $W$-invariant in $\l$, the equality (\ref{eq:keyequality-b-uno}) can be rewritten as
\begin{equation}
\label{eq:keyequality-b-due}
\int_{i\mathfrak a^*} a(\l)b(\l)\varphi_\l(x) \; \frac{d\l}{c(\l)c(-\l)} = \sum_{\mu \in \Lambda^+} (-1)^{|\mu|} d(\mu) a(\mu+\rho) \varphi_{\mu+\rho}(x)
\end{equation}
for all $x \in X_G$ sufficiently close to the base point $o=eK$.

Under the identifications $\l\equiv (\l_1,\dots,\l_l)$ by means of the basis $\Pi_*=\{\omega_1,\dots,\omega_l\}$ of $\mathfrak a^* \equiv \R^l$, the right-hand side of (\ref{eq:keyequality-b-due}) becomes
\begin{equation}
\label{eq:multipleseries}
\sum_{\mu_1=0}^{+\infty} \cdots \sum_{\mu_l=0}^{+\infty} (-1)^{|\mu|} a(\mu+\rho) d(\mu) \varphi_{\mu+\rho}(x)
\end{equation}
where
$$
\mu=(\mu_1,\dots,\mu_l)\,, \qquad \mu+\rho=(\mu_1+\rho_1,\dots,\mu_l+\rho_l)\,, \qquad
|\mu|=\mu_1+\dots+\mu_l\,.
$$
Moreover, the integral on the left-hand side becomes
\begin{equation} \label{eq:multipleintegral}
\int_{i\R} \cdots \int_{i\R} a(\l) \varphi_\l(x) \left( \frac{b(\l)}{c(\l)c(-\l)} \right) \; d\l_1\cdots d\l_l\,.
\end{equation}
Under suitable decay and convergence conditions, subsequent applications of the $1$-dimensional residue theorem to (\ref{eq:multipleintegral}) yields (\ref{eq:multipleseries}) provided:
\begin{enumerate}
\item
The function $\frac{b(\l)}{c(\l)c(-\l)}$ is meromorphic, with simple poles in the region $(\R^+)^l$
along the hyperplanes $\l_j=\mu_j+\rho_j$ with $\mu_j\in \Z^+$ and $j=1,\dots,l$.
\item
For $\mu=(\mu_1,\dots,\mu_l)\in (\Z^+)^l$, we have
$$(-2\pi i)^l \Res{\l_1=\mu_1+\rho_1} \dots \Res{\l_l=\mu_l+\rho_l}\, \frac{b(\l)}{c(\l)c(-\l)}=(-1)^{|\mu|} d(\mu)\,.$$
\end{enumerate}
Based on the rank-one case in \cite{Bertram}, we are therefore led to define $b(\l)$ by means of the equality
\begin{equation}
\label{eq:defb}
\frac{b(\l)}{c(\l)c(-\l)}= C_b d(\lambda-\rho) \prod_{j=1}^l \frac{1}{\sin\big(\pi(\l_j-\rho_j)\big)}
\end{equation}
where $C_b$ is a suitable constant. From the above arguments, we can compute that
\begin{equation}\label{eq:Cb}
C_b=\left(\frac{i}{2}\right)^l\,.
\end{equation}
To make the definition (\ref{eq:defb}) explicit, we need to further analyze the relation between the Plancherel density $[c(\l)c(-\l)]^{-1}$ for $X_G$ and the Plancherel density $d(\mu)$ for $X_U$.

\section{The Plancherel densities on $X_G$ and $X_U$}
\label{section:PlancherelDensities}
\noindent
Recall that if both $\b/2$ and $\b$ are roots, then $m_{\b/2}$ is even and $m_\b$ is odd. See e.g. \cite{He1}, Chapter X, Ex. F. 4. For a fixed $\b \in \Sigma_*^+$, the singularities of the function $[c_\b(\l)c_\b(-\l)]^{-1}$ are then described by distinguishing the following four cases:
\begin{enumerate}
\thmlist
\item $m_\b$ even, $m_{\b/2}=0$;
\item $m_\b$ odd, $m_{\b/2}=0$;
\item $m_\b$ odd, $m_{\b/2}/2$ even;
\item $m_\b$ odd, $m_{\b/2}/2$ odd.
\end{enumerate}
Recall the constant $\wrhob$ attached to $\beta \in \Sigma_*^+$ from formula (\ref{eq:wrhob}).
Following the computations yielding to formula (25) in \cite{HP},
we obtain the following lemma. See also \cite{Bertram}, Proposition 1.4.1.

\begin{Lemma}
\label{lemma:densitybeta}
Let $\b \in \Sigma_*^+$. Then
\begin{equation}\label{eq:dec-c_b}
\frac{1}{c_\b(\l)c_\b(-\l)}=C_\b\, p_\b(\l) q_\b(\l)
\end{equation}
where
\begin{enumerate}
\item
$C_\b$ is a positive constant (depending on $\b$ and on the multiplicities), explicitly given by
\begin{equation}\label{eq:Cbeta}
\text{$C_\b=4\pi \varepsilon(\b)$ \qquad where \qquad $\varepsilon(\beta)=\begin{cases}
(-1)^{m_\b/2} &\text{$m_\b$ is even}\\
(-1)^{(m_{\b/2}+m_\b-1)/2} &\text{$m_\b$ is odd}\end{cases}$}
\end{equation}
\item
$p_\b$ is a polynomial.
If $\wrhob>1/2$, then
\begin{align*}
p_\b(\l)&= \l_\b \big(\l_\b+\wrhob-1\big)\big(\l_\b+\wrhob-2\big) \cdots \big(\l_\b-(\wrhob-2)\big)
\big(\l_\b-(\wrhob-1)\big) \times \\
&\times \big(\l_\b+\big(\frac{m_{\b/2}}{4}-\frac{1}{2}\big)\big)\big(\l_\b+\big(\frac{m_{\b/2}}{4}-\frac{3}{2}\big)\big) \dots \big(\l_\b-\big(\frac{m_{\b/2}}{4}-\frac{3}{2}\big)\big)
\big(\l_\b-\big(\frac{m_{\b/2}}{4}-\frac{1}{2}\big)\big)\,,
\end{align*}
and the product on the second line does not occur if $m_{\b/2}=0$. If $\wrhob=1/2$, then
$p_\b(\l)=\l_\b$.
\item
$q_\b(\l)=1$ if $m_\b$ is even; if $m_\b$ is odd, then
\begin{equation}
q_\b(\l)=-\tan\big(\pi(\l_\b- \frac{m_{\b/2}}{4})\big)=\cot\big(\pi(\l_\b-\wrhob)\big)\,.
\end{equation}
\end{enumerate}
\end{Lemma}

The relation between the Plancherel measures on $X_G$ and $X_U$ is given by the following lemma.

\begin{Lemma}
\label{lemma:d-c}
The dimension $d(\mu)$ of the finite-dimensional
spherical representation of highest restricted weight $\mu \in \Lambda^+$ is given
by
$$d(\mu)=\left.\frac{c(\l-\mu)c(-\l+\mu)}{c(\l)c(-\l)}\right|_{\l=\mu+\rho}\,.$$
\end{Lemma}
\begin{proof}
This is Theorem 9.10, p. 321, in \cite{He3}.
\end{proof}

The apparent singularities in the formula in Lemma \ref{lemma:d-c} can be removed using Lemma
\ref{lemma:densitybeta} and the fact that the cotangent function is $\pi$-periodic. We obtain the following
formula from \cite{HP}, Proposition 3.5.

\begin{Prop}
\label{prop:d}
For $\mu\in \Lambda^+$ we have
$$d(\mu)=\frac{P(\mu+\rho)}{P(\rho)}\,$$
where
\begin{equation}
\label{eq:P}
P(\l)=\prod_{\b\in \Sigma_*^+} p_\b(\l)
\end{equation}
and $p_\b$ is the polynomial from Lemma \ref{lemma:densitybeta}.
\end{Prop}

According to Proposition \ref{prop:d}, the polynomial function on $\mathfrak a_\C^*$ extending $d(\mu)$ 
by means of Weyl integration formula can be written in terms of the positive unmultipliable roots as
\begin{equation}\label{eq:d-poly}
d(\l)=\frac{P(\l+\rho)}{P(\rho)}\,,  \qquad \l \in \mathfrak a_\C^*\,,
\end{equation}
where $P(\l)$ is as in (\ref{eq:P})\,.
\begin{Cor}
\label{cor:b}
\begin{enumerate}
\thmlist
\item
Let $C_b=\big(\frac{i}{2}\big)^l$ be the constant introduced in (\ref{eq:Cb}).
Then the function
\begin{equation}
\label{eq:bPlancherel}
\frac{b(\l)}{c(\l)c(-\l)}=\frac{C_b}{P(\rho)} \, P(\l) \; \prod_{j=1}^l \frac{1}{\sin\big(\pi(\l_j-\rho_j)\big)}
\end{equation}
is meromorphic on $\mathfrak a_\C^*$ with simple poles located along the hyperplanes
of equation
$$\pm \l_j -\rho_j=k_j$$
 where $k_j \in \Z^+$ and $j\in\{1,\dots,l\}$.
\item
For all $w \in W$ the function $\frac{b(w\l)}{c(\l)c(-\l)}$ is holomorphic on the tube
$T_1={\rm C}(\rho)^0+i \mathfrak a^*$ of Lemma \ref{lemma:Trho}.
\item
For all $\l \in \mathfrak a_\C^*$ we have
\begin{equation}
\label{eq:bT}
b(\l)=K_b T(\l) \; \prod_{j=1}^l \frac{1}{\sin\big(\pi(\l_j-\rho_j)\big)}
\end{equation}
where
\begin{equation}
\label{eq:T}
T(\l)=\prod_{\b\in \Sigma_*^+} t_\b(\l)
\end{equation}
and
\begin{equation}
\label{eq:tbeta}
t_\b(\l)=\begin{cases} 1 &\text{if $m_\b$ is even}\\
                      \tan\big(\pi(\l_\b-\wrhob)\big) &\text{if $m_\b$ is odd}
\end{cases}
\end{equation}
is the inverse of the function $q_\b$ from Lemma \ref{lemma:densitybeta}. Moreover,
$K_b$ is a constant depending on the multiplicities. It is given explicitly by
$$K_b=\left(\frac{i}{2}\right)^l \frac{c_0^2}{P(\rho)} \, \prod_{\b \in \Sigma_*^+} C_\b\,,$$
where $C_\b$ is as in (\ref{eq:Cbeta}) and $c_0$ is the constant appearing in the definition (\ref{eq:c}) of
Harish-Chandra's $c$-function.
\end{enumerate}
\end{Cor}
\begin{proof}
Immediate consequence of (\ref{eq:defb}), Lemma \ref{lemma:densitybeta} and Proposition \ref{prop:d}.
For the list of singular hyperplanes, notice that for fixed $j\in \{1,\dots,l\}$, the zeros of the function $\sin\big(\pi(\l_j-\rho_j)\big)$ are located along the hyperplanes of equation
$\l_j-\rho_j=k_j$ with $k_j \in \Z$. Recall that $\l_j=\l_{\b_j}$ and $\widetilde{\rho}_{\b_j}=\rho_j$.
The polynomial $p_{\b_j}$ is divided by
$$\big(\l_j-(\rho_j-1)\big)\big(\l_j-(\rho_j-2)\big)\cdots \big(\l_j+(\rho_j-2)\big)\big(\l_j+(\rho_j-1)\big)\,.$$
Hence all singularities of $[\sin\big(\pi(\l_j-\rho_j)\big)]^{-1}$ for $|\l_j|<\rho_j$ are canceled.
No other singularities of this function are canceled by zeros of $p_{\b_j}$.

Because of (a), the function $\frac{b(\l)}{c(\l)c(-\l)}$ is holomorphic on the tube
$T'_1$ of Lemma \ref{lemma:Trho}. Hence $\frac{b(w\l)}{c(\l)c(-\l)}$ is holomorphic
on the largest $W$-invariant tube domain contained in $T'_1$, which is $T_1=\cap_{w \in W} w(T'_1)$.

The formula for the constant $K_b$ is obtained by comparing (\ref{eq:bPlancherel}), (\ref{eq:bT}), (\ref{eq:c}) and
(\ref{eq:dec-c_b}).
\end{proof}

\begin{Rem}
\label{rem:b}
According to the four cases for $m_\b$ and $m_{\b/2}$ listed in Section \ref{section:PlancherelDensities}, we have $\wt\rho_\beta \in \Z$
in cases (a) and (d), and $\wt\rho_\beta \in \Z+1/2$ in cases (b) and (c).
We can therefore write (\ref{eq:bT}) as
\begin{equation}
\label{eq:bexplicit}
b(\l)=K'_b \Big(\prod_{\stackrel{\b \in\Sigma_*^+\setminus\{\b_1,\dots,\b_l\}}{\text{cases (b) or (c)}}} \hskip -6mm  \cot(\pi\l_\b)\Big)
\Big( \prod_{\stackrel{\b \in\Sigma_*^+\setminus\{\b_1,\dots,\b_l\}}{\text{case (d)}}}  \hskip -6mm 
\tan(\pi\l_\b) \Big) 
\Big( 
\hskip -4mm
\prod_{\stackrel{j\in\{1,\dots,l\}}{\text{cases (a),(b) or (c)}}}  \hskip -4mm   \frac{1}{\sin(\pi\l_j)}\Big)
\Big( \prod_{\stackrel{j\in\{1,\dots,l\}}{\text{case (d)}}} \!   \frac{1}{\cos(\pi\l_j)} \Big)
\end{equation}
where $K'_b=\pm K_b$ and the sign depends on the parity of the multiplicities.
\end{Rem}

Remark \ref{rem:b} immediately implies the following corollary. Notice that the cases (a), (b) and (c) for the root multiplicities correspond to the situation in which $(m_{\beta/2})/2$ is even.

\begin{Cor}
\label{cor:Pib}
Let
\begin{equation}
\label{eq:Pi}
\Pi(\l)=\prod_{\b\in\Sigma_*^+} \l_\b\,
\end{equation}
and let $T_{\Sigma,m}$ be as in (\ref{eq:TSigmam}).
Then $\Pi(\l)b(\l)$ is holomorphic on $T_{\Sigma,m}$.
\end{Cor}

\begin{Ex}[The rank-one case]
In the (real) rank-one case, $\mathfrak a$ is one dimensional. Then $\Sigma^+$ consists of at
most two elements: $\beta$ and, possibly, $\beta/2$. Hence $\Sigma_*^+=\{\beta\}$, $l=1$ and $\beta_1=\beta$. According to Remark \ref{rem:b}, we have
$$b(\l)=\begin{cases} \pm K_b [\sin(\pi\l_1)]^{-1} &\text{in cases (a),(b) and (c)}\\
         \pm K_b [\cos(\pi\l_1)]^{-1} &\text{in case (d)}.
        \end{cases}.$$
This case has been previously considered by Bertram in \cite{Bertram}.
\end{Ex}

\begin{Ex}[The even multiplicity case]
Suppose that $\Sigma$ is a reduced root system and that all roots multiplicities are even. Geometrically,
even multiplicities correspond to Riemannian symmetric spaces of the noncompact type $G/K$ with the property that all Cartan subalgebras in the Lie algebra $\mathfrak g$ of $G$ are conjugate under the adjoint group of $\mathfrak g$. The simplest examples occur when $\mathfrak g$ admits a complex structure, in which case all root multiplicities are equal to $2$.  In the even multiplicity case, one has $b(\l)=K_b \prod_{j=1}^l [\sin(\pi \l_j)]^{-1}$. In the complex case, one can compute that $K_b=(i/2)^l$. In particular, for the complex rank-one case corresponding to the pair of symmetric spaces $G/K=\SL(2,\C)/\SU(2)$ and $U/K\cong \SU(2)$, we have $b(\l)=\frac{i}{2} \, [\sin(\pi \l_1)]^{-1}$. Thus, in this case and with our choice of the coordinate $\l_1$ in $\mathfrak a_\C^*$,
the function $b$ agrees with the one of Ramanujan for the case of $\R^+$.
\end{Ex}

\begin{Rem}
In comparison to the classical version of the Master Theorem, the assumption $P>0$ in the statement of Ramanujan's Theorem for symmetric spaces strongly restricts the class of spherical Fourier series to which our theorem can be applied. On the other hand, by replacing the pair $({\rm U}(1),\R^+)$ by a pair of Riemannian symmetric spaces in duality, we are considering a much richer class of different geometric situations where our theorem applies. According to the choice of the dual pair of symmetric spaces, the restrictions to the Cartan subspace $A$ of the spherical functions $\varphi_{\mu+\rho}$ provide several classes of orthogonal polynomials of Jacobi type in several variables. Likewise, by the integral formulas corresponding to the decomposition $G=KAK$, the spherical Fourier transforms gives different specializations of Jacobi transforms in several variables. The simplest example of these specializations corresponds to the complex case considered above. In this case, $G=K_\C$ and $U=K\times K$ where $K$ is a compact connected semisimple Lie group. The corresponding Riemannian symmetric spaces in duality are $X_G=K_\C/K$ and $X_U=(K\times K)/K$ . The space $X_U$ can be identified with $K$. In this way, the $K$-invariant functions on $X_U$ correspond to the central functions on $K$. To simplify notation, we assume in the following that $K$ is simply connected. The spherical representations of $U$ are of the form $\pi_\mu=\delta_\mu \otimes \overline{\delta_\mu}$ where $\mu \in \Lambda^+$ is the highest weight of the irreducible representation $\delta_\mu$ of $K$ and $\overline{\delta_\mu}$ denotes the contregradient representation of $\delta_\mu$. The spherical functions on $X_U$ are therefore the normalized characters $\frac{1}{\dim \delta} \chi_\mu$ where $\chi_\mu$ is the character of $\delta_\mu$. On the noncompact side, one can obtain explicit formulas by using the integral formulas for the $KAK$ decomposition of $G$ and the explicit formulas on $A$ for the
spherical functions on Riemannian symmetric spaces of the noncompact type and $G$ with complex structure; see e.g. \cite{He2}, Ch. I, Theorem 5.8 and Ch. IV, Theorem 5.7. Since $\Sigma=\Sigma_*$ is reduced and all multiplicities are equal to $2$, one has $\rho=\sum_{\b \in \Sigma^+} \b$.
Set $$f(x)=\sum_{\mu \in \Lambda^+} (-1)^{|\mu|} a(\mu +\rho) \chi_\mu(x)\, \qquad x \in K_\C$$
with $a \in \mathcal H(A,P,\delta)$.
Then Ramanujan's interpolation formula (\ref{eq:Ramanujansymm}) becomes the Fourier integral formula
$$\int_{\frak a^+} f(\exp H) \Delta(H) e^{\l(H)} \; dH =\sum_{w \in W} a(w\l)b(w\l)$$
where
\begin{align*}
\Delta(H)&=\prod_{\b \in \Sigma^+} (e^{\b(H)}-e^{-\b(H)})= \sum_{w \in W} (\det w) e^{w \rho(H)}\,,\\
b(\l)&=c_b \prod_{j=1}^l \sin(\pi \l_j)^{-1}
 \end{align*}
and $c_b$ is a suitable normalizing constant.
\end{Rem}

\section{Some estimates}
\label{section:estimates}
\noindent
In this section we collect some estimates which will be needed in the proof of Theorem \ref{thm:RamanujanSymm}.

For every $\l\in\mathfrak a_\C^*$ the spherical function $\varphi_\l$ extends holomorphically as a $K_\C$-invariant function
on the domain $K_\C \exp(2\Omega_\pi)\cdot o$ in $X_\C$, where
\begin{equation}
\label{eq:Omegapi}
\Omega_\pi=\{H\in \mathfrak a_\C: \text{$|\b(\Im H)|<\pi/2$ for all $\b\in\Sigma$}\}\,.
\end{equation}
The estimates of the holomorphically extended spherical functions that are given in Lemma \ref{lemma:estphil} below will be sufficient to our purposes.

Recall that $\mathfrak a^*_+=\{\l\in \mathfrak a^* : \text{$\l_\b\geq 0$ for all $\b\in\Sigma_+^*$}\}$.
Notice that $\mathcal H(\delta) \supset \mathfrak a^*_+ + i\mathfrak a^*$ for all $\delta>0$.
Recall also the constant $\Omega$ from (\ref{eq:Omega}).

\begin{Lemma}
\label{lemma:estphil}
There is a constant $C>0$ so that
\begin{equation}
\label{eq:OpdamEstimates}
|\varphi_\l(\exp H \cdot o)|\leq C  e^{-\min_{w \in W} \Im(w\l(H_2))+\max_{w\in W} \Re(w\l(H_1))}
\end{equation}
for all $\l \in \mathfrak a_\C^*$ and all
$H=H_1+iH_2 \in \overline{\Omega_\pi}$ with $H_1, H_2 \in \mathfrak a$.
In particular:
\begin{enumerate}
\thmlist
\item
for all $\l \in \mathfrak a^*_+ + i\mathfrak a^*$ and $H \in \mathfrak a$ we have
$$|\varphi_\l(\exp H \cdot o)|\leq C  e^{\Omega\|H\|(\sum_{j=1}^l \Re\l_j)}$$
\item
for all  $H\in \overline{\Omega_\pi}$ and $\l\in \mathfrak a^*$ we have
$$|\varphi_\l(\exp H \cdot o)|\leq C  e^{\|\Im H\|\|\Im \l\|}\,.$$
\end{enumerate}
\end{Lemma}
\begin{proof}
Estimates (\ref{eq:OpdamEstimates}) are due to Opdam; see \cite{OpdamActa}, Proposition 6.1(2) and Theorem 3.15. For (a), we can suppose by $W$-invariance that $H \in \overline{\mathfrak a^+}$.
In this case, for $\l \in \mathfrak a^*_+ + i\mathfrak a^*$, we have
\begin{equation*}
0\leq \Re\l(H)=\sum_{j=1}^l \Re\l_j \omega_j(H) \leq \Omega \|H\| \big(\sum_{j=1}^l \Re\l_j\big)\,.
\end{equation*}
Part (b) follows immediately from (\ref{eq:OpdamEstimates}).
\end{proof}

\begin{Lemma}
\label{lemma:estP}
Let $P$ be as in (\ref{eq:P}). Then there are positive constants $C_0, C'_0$ and $C''_0$ so that
\begin{equation*}
|P(\l)| \leq C_0 \prod_{\b\in\Sigma^+_*} (1+|\l_\b|)^{m_{\b/2}+m_\b}
\leq C'_0 (1+\|\l\|)^M \leq C_0'' \prod_{j=1}^l (1+|\l_j|)^M
\end{equation*}
where
\begin{equation}
\label{eq:M}
M=\sum_{\b\in\Sigma_*^+} (m_{\b/2}+m_\b)\,.
\end{equation}
\end{Lemma}
\begin{proof}
The first inequality is an immediate consequence of the formula for $p_\b$ in Lemma \ref{lemma:densitybeta}, which gives $p_\b$ as polynomial of degree $m_{\b/2}+m_\b$ in $\l_\b$. For the second, notice that $|\l_\b|\leq \|\b\|^{-1} \|\l\|$. The final inequality follows immediately from $\l=\sum_{j=1}^l \l_j \omega_j$.
\end{proof}

We define
\begin{equation}
\label{eq:Q}
Q(\l)=\prod_{\b\in\Sigma_*^+} q_\b(\l)
\end{equation}
where $q_\b$ is the function defined in Lemma \ref{lemma:densitybeta}.

Observe that there is a constant $K>0$ so that
\begin{equation}
\label{eq:est-sin}
\big|\sin\big(\pi(\l_j-\rho_j)\big)\big|^{-1} \leq K e^{-\pi|\Im \l_j|}
\end{equation}
for $|\Im \l_j|\geq 1$ or for $\Re\l_j=\rho_j+N+1/2$ with $N\in \Z^+$.

The following lemma contains the estimates needed to apply the Residue Theorem.

\begin{Lemma}
\label{lemma:estimates-residues}
\begin{enumerate}
\thmlist
\item
Let $N$ be a positive integer and let $M$ be as in (\ref{eq:M}).
Let $\l=\sum_{j=1}^l \l_j \omega_j \in\mathfrak a_\C^*$ with
$|\Im \l_j|\geq 1$ or $\Re\l_j=\rho_j+N+1/2$ or $\Re\l_j=0$ for all $j=1,\dots,l$.
Then there is a positive constant $C_1$, independent of $N$, so that
\begin{equation*}
\left|\frac{b(\l)}{c(\l)c(-\l)}\right| \leq C_1 \prod_{j=1}^l \left[(1+|\l_j|)^M e^{-\pi |\Im \l_j|}\right]\,.
\end{equation*}
\item
Set
\begin{multline}
\label{eq:B}
B=\Big\{\l=\sum_{j=1}^l \l_j \omega_j \in\mathfrak a_+^*+i\mathfrak a^*:
\text{$|\Im \l_j|\geq 1$ or}\\ \text{$\Re\l_j \in (\rho_j+\Z^+ +1/2) \cup \{0\}$ for all $j=1,\dots,l$}\Big\}.
\end{multline}
Let $a \in \mathcal H(A,P,\delta)$. Then there is a constant $C_2>0$ so that for all $\l\in B$ and $H \in \overline{\mathfrak a^+}$ we have
\begin{equation}
\label{eq:main-est-shift}
\left|\frac{a(\l)b(\l)}{c(\l)c(-\l)} \varphi_\l(\exp H\cdot o) \right|
\leq C_2 \prod_{j=1}^l \left[ (1+|\l_j|)^M e^{(A-\pi)|\Im \l_j| +(\|H\|\Omega-P) \Re\l_j} \right]\,.
\end{equation}
\end{enumerate}
\end{Lemma}
\begin{proof}
For $|\Im \l_j|\geq 1$ or $\Re\l_j=\rho_j+N+1/2$,
the estimate in (a) is a consequence of (\ref{eq:bPlancherel}), (\ref{eq:est-sin}) and Lemma \ref{lemma:estP}. The inequality holds also if $\Re\l_j=0$ for some $j$, as the possible singularity
of $\sin\big(\pi(\l_j-\rho_j)\big)^{-1}$ at $\l_j=0$ is cancelled by the factor $\l_j$ in
$p_{\b_j}(\l)$; see formulas (\ref{eq:bPlancherel}) and Lemma \ref{lemma:densitybeta}.
Part (b) follows from (a) and Lemma \ref{lemma:estphil},(a).
\end{proof}

The next lemma will be useful to prove the independence on $\sigma \in B(T_\delta)$ for the integral occurring in Part (b) of Ramanujan's Theorem.

\begin{Lemma}
\label{lemma:estimatesTdelta}
Let $0<\delta \leq 1$ and let $T_\delta$ be the tube domain from (\ref{eq:Trhouno}). Let $M$ be the
constant defined in (\ref{eq:M}).
\begin{enumerate}
\thmlist
\item
There is a constant $C_\delta >0$ so that
\begin{equation}
\label{eq:estimatebcTdelta}
\left|\frac{b(\l)}{c(\l)c(-\l)} \right| \leq C_\delta (1+\|\l\|)^M e^{-\pi\big(\sum_{j=1}^l
|\Im \l_j|\big)}
\end{equation}
for all $\l \in T_\delta$.
\item
Let $a \in \mathcal H(A,P,\delta)$. For every $R>0$ and every integer $N\geq 0$ there is a constant $C_{R,N,\delta} >0$ so that for all $\l\in T_\delta$ and $H \in \mathfrak a$ with $\|H\|< R$, we have
\begin{equation}
\label{eq:estabphilN}
\left|\frac{a(\l)b(\l)}{c(\l)c(-\l)} \varphi_\l(\exp H\cdot o) \right|
\leq
C_{R,N,\delta} (1+\|\l\|)^{-N}
\end{equation}
Consequently,
\begin{equation}
\label{eq:estWabphilN}
\left|\big(\sum_{w\in W} a(w\l)b(w\l)\big)\frac{\varphi_\l(\exp H\cdot o)}{c(\l)c(-\l)} \right|
\leq
C_{R,N,\delta} |W|(1+\|\l\|)^{-N}\,.
\end{equation}
\end{enumerate}
\end{Lemma}
\begin{proof}
The polynomial
$$p_j(\l)=\big(\l_j-(\rho_j-1)\big)\big(\l_j-(\rho_j-2)\big)\cdots \big(\l_j-(-\rho_j+2)\big)\big(\l_j-(-\rho_j+1)\big)$$ is a divisor of $p_{\b_j}(\l)$. Hence
\begin{align}
\frac{b(\l)}{c(\l)c(-\l)}&=\frac{C_b}{P(\delta)}\, P(\l) \prod_{j=1}^l \frac{1}{\sin(\pi(\l_j-\rho_j))} \notag\\
&=\frac{C_b}{P(\delta)}\, \wt{P}(\l) \prod_{j=1}^l \frac{p_j(\l)}{\sin(\pi(\l_j-\rho_j))}
 \label{eq:estTdelta}
\end{align}
for a certain polynomial $\wt P$. For fixed $\eta \in ]0,1[$, the function $\frac{z}{\sin(\pi z)}$
is bounded on $\{z \in \C: |\Im z|\leq 1, |\Re z|\leq \eta\}$. By (\ref{eq:est-sin}), we conclude that there is a constant $C'_\delta>0$ so that for any fixed $j=1,\dots, l$ and every $\l=\sum_h \l_h \omega_h$ with $|\Re\l_j|\leq \delta \rho_j$ and arbitrary $\l_h \in \C$ with $h\neq j$, we have
$$\left| \frac{p_j(\l)}{\sin\big(\pi(\l_j-\rho_j)\big)}\right| \leq C'\delta (1+|\l_j|)^{\deg p_j} e^{-\pi|\Im \l_j|}\,.$$
We obtain Part (a) from these estimates and (\ref{eq:estTdelta}).

To prove (b), observe first that by Lemma \ref{lemma:Trho}, the function
$\frac{b(w\l)}{c(\l)c(-\l)}$ is holomorphic on $T_1 \supset T_\delta$. Moreover, $T_\delta$ is a $W$-invariant subset of $\mathcal H(\delta)$. Hence $\frac{a(w\l)b(w\l)}{c(\l)c(-\l)} \varphi_\l(\exp H\cdot o)$ is holomorphic on $T_\delta$.
Notice that $|\varphi_\l(\exp H \cdot o)| \leq \varphi_{\Re\l}(\exp H\cdot o)$. Let $R>0$ be fixed. Since the basis $B(T_\delta)$ of the tube domain $T_\delta$ has compact closure, it follows by continuity, that there is a constant $C_{R,\delta}>0$ so that
\begin{equation}
\label{eq:estphilcomp}
|\varphi_\l(\exp H \cdot o)| \leq C_{R,\delta}
\end{equation}
for all $\l \in T_\delta$ and $H\in \mathfrak a$ with $\|H\| \leq R$.
(This can also be obtained from Lemma \ref{lemma:estphil}.)

Suppose that $a\in\mathcal H(A,P,\delta)$. By (\ref{eq:estphilcomp}) and Part (a), there is a constant
$C'_{R,\delta}>0$ so that
\begin{equation*}
\left|\frac{a(\l)b(\l)}{c(\l)c(-\l)} \varphi_\l(\exp H\cdot o) \right|
\leq
C'_{R,\delta} (1+\|\l\|)^{M} e^{(A-\pi)\sum_{j=1}^l |\Im \l_j|}\,.
\end{equation*}
This implies (\ref{eq:estabphilN}) as $A < \pi$.
\end{proof}

We shall also need the following result, which is a local version of a classical argument by Malgrange (see
\cite{He3}, p. 278, and \cite{P}, Lemma 4.2).

\begin{Lemma}
\label{lemma:localMalgrange}
Let $V^\prime $ be an open domain in $\mathfrak a_\C^*$ and let $H:V^\prime \to \C$ a holomorphic function satisfying the following property: there exist constants $R \in \R$, $s>0$ and $C>0$ so that
$$|H(\l)|\leq C(1+\|\l\|)^s e^{R\|\Im \l\|} $$
for all $\l \in V^\prime $.
Let $\tau >0$ and let $V$ be an open domain in  $\mathfrak a_\C^*$ such that
$$V_\tau:=\big\{\nu \in \mathfrak a_\C^*: \text{$\exists \l \in V$ with $\|\l-\nu\|\leq \tau$}\big\} \subset V^\prime \,.$$
Let $p$ be a polynomial such that $F(\l)=\frac{H(\l)}{p(\l)}$ is holomorphic on $V^\prime $.
Then there is a constant $C_\tau>0$ (depending also on $C,R,s$) such that
$$|F(\l)|\leq C_\tau (1+\|\l\|)^s e^{R\|\Im \l\|} $$
for all $\l \in V$.
\end{Lemma}
\begin{proof}
Let $m=\deg p$. By Cauchy's integral formula, for any multiindex $\alpha$ there is a constant
$C_{m,\a}>0$ so that for every $\l \in V$ we have
$$|F(\l)(\partial^\alpha p)(\l)| \leq C_{m,\alpha} \int_{\|\xi\| \leq \tau} |F(\l+\xi)||p(\l+\xi)| \; d\xi\,.$$
Choose $\alpha$ so that $\partial^\alpha p$ is a constant $d\neq 0$. We obtain:
\begin{align*}
|F(\l)|& \leq d^{-1} C_{m,\alpha} \int_{\|\xi\| \leq \tau} |H(\l+\xi)| \; d\xi\\
& \leq d^{-1} C C_{m,\alpha} \int_{\|\xi\| \leq \tau} (1+\|\l+\xi\|)^s e^{R(\|\Im \l+\Im \xi\|)} \; d\xi\\
& \leq C_{\tau} (1+\|\l\|)^s e^{R\|\Im \l\|}\,,
\end{align*}
where $C_\tau=d^{-1} C C_{m,\alpha} (1+\tau)^s e^{|R|\tau} \int_{\|\xi\| \leq \tau} \,d\xi$.
\end{proof}

The following lemma shows that if $a\in \mathcal H(A,P,\delta)$, then there is a constant
$\varepsilon \in ]0,1[$ (depending on $A,P,\delta$) so that the function
\begin{equation}
\label{eq:atilde}
\wt a(\l)=\sum_{w\in W} a(w\l) b(w\l)
\end{equation}
belongs to the $W$-invariant Schwartz space $\mathcal S(\mathfrak a_\varepsilon^*)^W$ on the tube domain $T_\varepsilon$ around $i\mathfrak a^*$; see (\ref{eq:Trho}) for the definition of $T_\varepsilon$ and
section \ref{subsection:spherical} for the definition of the $W$-invariant Schwartz space.

For $0 \leq \eta <1/2$ set
\begin{multline*}
T_{\Sigma,m,\eta}=\{ \l \in \mathfrak a_\C: \text{$|\Re\l_\b|<1-\eta$, $\b\in \Sigma_*^+$ with $(m_{\b/2})/2$ even}\} \\
\cap  \{ \l \in \mathfrak a_\C: \text{$|\Re\l_\b|<1-\eta$, $\b\in \Sigma_*^+$ with $(m_{\b/2})/2$ odd}\}\,.
\end{multline*}
So $T_{\Sigma,m,0}=T_{\Sigma,m}$ is the tube domain on which $\Pi(\l)b(\l)$ is holomorphic; see Corollary \ref{cor:Pib}.

\begin{Lemma}
\label{lemma:Schwartztube}
\begin{enumerate}
\thmlist
Set $s=|\Sigma_*^+|$.
\item
Let $0<\eta<1/2$. Then there is a constant $C_\eta>0$ so that
$$|\Pi(\l)b(\l)|\leq C_\eta (1+\|\l\|)^s e^{-\pi\big(\sum_{j=1}^l |\Im \l_j|\big)}$$
for all $\l \in T_{\Sigma,m,\eta}$.
\item
Let $a \in \mathcal H(A,P,\delta)$ and set $\wt a(\l)=\sum_{w \in W} a(w\l)b(w\l)$. Then $\wt a$ is holomorphic in
$T_{\Sigma,m} \cap T_\delta$. Moreover, let $0 <\eta <\min\{1/2,\delta\}$. Then there is a constant $C_{\eta,a}>0$ so that
$$|\wt a(\l)| \leq C_{\eta,a} (1+\|\l\|)^s e^{(A-\pi)c_2 \|\Im \l\|}$$
for all $\l\in T_{\Sigma,m,\eta} \cap T_{\delta-\eta}$. Here $c_2$ is the positive constant introduced in
Remark \ref{rem:equivalentnorms}.
\item
Let $$\gamma=\min\Big\{\delta, \min_{b\in\Sigma_*^+, \text{$(m_{\b/2})/2$ even}} (\wt\rho_\b)^{-1},
\min_{\b\in\Sigma_*^+, \text{$(m_{\b/2})/2$ odd}} (2\wt\rho_\b)^{-1}\Big\}\,\in ]0,1[\,.$$
Then $T_\gamma \subset T_{\Sigma,m,\eta} \cap T_\delta$.
Moreover, let $0<\varepsilon <\gamma$. Then $\wt a \in \mathcal S(\mathfrak a^*_\varepsilon)^W$, the $W$-invariant Schwartz space on
the tube domain $T_\varepsilon$.
\end{enumerate}
\end{Lemma}
\begin{proof}
Since $\Pi(\l)b(\l)$ is bounded on $T_{\Sigma,m,\eta}$, the proof of the estimate in (a) follows the same argument as in part (a) of Lemma \ref{lemma:estimatesTdelta}.

To prove part (b), notice first that, by Corollary \ref{cor:Pib}, on $T_{\Sigma,m}$
the function $b(\l)$ has at most simple poles on hyperplanes of the form $\l_\b=0$ with $\b \in\Sigma_*^+$. The same property holds on $T_{\Sigma,m} \cap T_\delta$ for $b(w\l)a(w\l)$, with $w \in W$, and hence for $\wt a(\l)$. But $\wt a(\l)$, as $W$-invariant function, cannot admit first order singularities on root hyperplanes through the origin. Thus $\wt a$ is holomorphic on $T_{\Sigma,m} \cap T_\delta$.

Let $0<\eta<\eta'<\min\{1/2,\delta\}$. By (a), there is a constant $C_{\eta'}>0$ so that
$$|\Pi(\l)b(\l)a(\l)| \leq C_{\eta'} (1+\|\l\|)^s e^{(A-\pi)\|\Im \l\|_1 }
\leq C_{\eta'}  (1+\|\l\|)^s e^{(A-\pi)c_2\|\Im \l\|}$$
for all $\l\in T_{\Sigma,m,\eta'}$.
Since $\Pi(\l)$ is $W$-skew-invariant, we obtain, on $T_{\Sigma,m,\eta'}$:
$$|\pi(\l)\wt a(\l)| \leq C_{\eta'} |W| (1+\|\l\|)^s e^{(A-\pi)c_2 \|\Im \l\|}\,.$$
The estimate for $\wt a$ on $T_{\Sigma,m,\eta}$ follows then from Lemma  \ref{lemma:localMalgrange}.

To show that $T_\gamma \subset T_{\Sigma,m} \cap T_\delta$, notice that $\gamma \leq \delta$ and that if
$\l\in T_{\gamma}$, then $|\Re\l_\b|<\gamma \wt \rho_\b \leq 1$ if $(m_{\b/2})/2$ is even and $\leq
1/2$ if $(m_{\b/2})/2$ odd.
Hence $T_\gamma \subset T_{\Sigma,m}$.
The property that $\wt a \in \mathcal S(\mathfrak a^*_\varepsilon)^W$ for $0<\varepsilon < \gamma$ is a consequence of Cauchy's estimates. Indeed, let $\wt\varepsilon=(\gamma-\varepsilon)/2$. For
$\l\in T_\varepsilon$, let $D=\{\nu \in \mathfrak a^*_\C: \text{$|\nu_j-\l_j|\leq \wt\varepsilon$
for all $j=1,\dots, l$}\}$ be the closed polydisc with center $\l=(\l_1,\dots,\l_l)$ and multiradius $(\varepsilon,\dots, \varepsilon)$. Then $D \subset T_\gamma$. According to Cauchy's estimates (see e.g. \cite{Krantz}, Lemma 2.3.9), for every multiindex $\alpha$,
$$|\partial^\a \wt a(\l)| \leq \frac{\a!}{\wt \varepsilon^{\,|\a|}} \, \sup_{\nu \in D} |\wt a(\nu)|\,.$$
Since
$|\wt a(\nu)|\leq C_\gamma (1+\|\nu\|)^s e^{(A-\pi)c_2\|\Im \nu\|}$, by estimating
$\|\nu\|$ and $\|\Im \nu\|$ in terms of $\|\l\|$ and $\|\Im \l\|$, respectively, we obtain for a constant $C_{\wt\varepsilon, \gamma}$ depending on $\wt\varepsilon$ and $\gamma$ but not on $\lambda$:
$$|\wt a(\l)|\leq C_{\wt\varepsilon,\gamma} (1+\|\l\|)^s e^{(A-\pi)c_2\|\Im \l\|}\,.$$
Since $A<\pi$, we conclude the required rapid decay.
\end{proof}

\section{Proof of Ramanujan's Master theorem for symmetric spaces}
\label{section:proofRamanujanSymm}
\noindent
In this section we prove Theorem \ref{thm:RamanujanSymm}.
Part 1  is an immediate consequence of Lassalle's Theorem \ref{thm:Lassalle} with $F(\mu)=(-1)^{|\mu|} a(\mu+\rho)$ and $\varepsilon=P/\Omega$. Indeed, for $\mu \in \Lambda^+$ we have
$\Omega \big(\sum_{j=1}^l \mu_j\big) = \Omega \|\mu\|_1 \geq \|\mu\|\,.$
Hence, for $a\in\mathcal H(\delta)$,
$$|a(\mu+\rho)|\leq C \prod_{j=1}^l e^{-P(\mu_j+\rho_j)}=C'e^{-P \sum_{j=1}^l \mu_j} \leq C' e^{-\frac{P}{\Omega}\|\mu\|}\,,$$
where $C'=C \prod_{j=1}^le^{-P\rho_j}$.

To prove Part 2, let $N$ be a positive integer. For $j=1,\dots,l$, let $C_{j,N}$ be the closed rectangular contour in
the $\l_j$-plane passing clockwise through its vertices $-iN$, $iN$, $iN+\big(\rho_j+N+1/2\big)$ and
$-iN+\big(\rho_j+N+1/2\big)$.
Set
$$f(\l)=C_b a(\l) \varphi_\l(x)d(\l-\rho)\,.$$
Recall that the spherical function $\varphi_\l$ is an entire function of $\l \in \mathfrak a_\C^*$.
Hence $f$ is holomorphic on $\mathcal H(\delta)$.
Suppose $\l_2,\dots,\l_l$ are fixed values so that $\l_j-\rho_j\notin \Z^+$.
Then, by the residue theorem, we have
\[\oint_{C_{1,N}} f(\l) \prod_{j=1}^l \frac{1}{\sin\big(\pi(\l_j-\rho_j)\big)} \; d\l_1=
(-2\pi i) \sum_{\mu_1=0}^N \Res{\l_1=\rho_1+\mu_1} \left( f(\l) \prod_{j=1}^l \frac{1}{\sin\big(\pi(\l_j-\rho_j)\big)}\right) \]
\[\hbox to 9em{}=(-2\pi i) \sum_{\mu_1=0}^N f(\rho_1+\mu_1, \l_2,\dots,\l_l) \frac{(-1)^{\mu_1}}{\pi} \,
\prod_{j=2}^l \frac{1}{\sin\big(\pi(\l_j-\rho_j)\big)}\,.
\]
Iterating, we obtain
\begin{align*}
\oint_{C_{1,N}}\cdots \oint_{C_{l,N}} f(\l) &\prod_{j=1}^l \frac{1}{\sin\big(\pi(\l_j-\rho_j)\big)} \; d\l_1 \cdots d\l_l =\\
&=(-2\pi i)^l \sum_{\mu_1=0}^N \cdots \sum_{\mu_l=0}^N f(\rho_1+\mu_1, \dots,\rho_l+\mu_l) \frac{(-1)^{\mu_l+\cdots+\mu_l}}{\pi^l}\\
&=(-2i)^l C_b \sum_{\mu_1=0}^N \cdots \sum_{\mu_l=0}^N (-1)^{|\mu|} d(\mu)a(\rho+\mu)
\varphi_{\rho+\mu}(x).
\end{align*}
Thus, by (\ref{eq:Cb}),
\begin{equation}\label{eq:residuesN}
\oint_{C_{1,N}}\cdots \oint_{C_{l,N}} a(\l)\varphi_\l(x) \frac{b(\l)}{c(\l)c(-\l)} \; d\l_1 \cdots d\l_l
=\sum_{\mu_1=0}^N \cdots \sum_{\mu_l=0}^N (-1)^{|\mu|} d(\mu)a(\rho+\mu)
\varphi_{\rho+\mu}(x)\,.
\end{equation}

By Part 1, the right-hand side of (\ref{eq:residuesN}) converges as $N \to \infty$ to
$\sum_{\mu\in \Lambda^+} (-1)^{|\mu|} d(\mu) a(\mu+\rho) \varphi_{\mu+\rho}(x)$ for
$x=\exp H\cdot o$ with $\|H\|<P/\Omega$, the convergence being normal on compacta of $\exp((P/\Omega)B)\cdot o$.

For the limit of the left-hand side of (\ref{eq:residuesN}), we shall use the estimate (\ref{eq:main-est-shift}). Observe first that the domain of integration for the left-hand side of (\ref{eq:residuesN})
is
$$D(C_{1,N}, \dots, C_{2,N})=\big\{\l=\sum_{j=1}^l \l_j \omega_j \in \mathfrak a^*_\C: \text{$\l_j \in C_{j,N}$ for all $j=1,\dots,l$}\big\}$$
and
$$\bigcup_{N=1}^\infty D(C_{1,N}, \dots, C_{2,N}) \subset B$$
where $B$ is the set on which the estimate (\ref{eq:main-est-shift}) holds.

\begin{Lemma}
\label{lemma:integrals}
Let $j\in\{1,\dots,l\}$ be fixed. For a positive integer $N$, let
$\gamma_{j,N}$ be the portion of the contour $C_{j,N}$ from $iN$ to $-iN$, and let
$\eta_{j,N}$ be the vertical portion of $C_{j,N}$ from $-iN$ to $iN$
(see Figure 1).
Let $\tau<0$ and $\sigma <0$ be fixed constants.
Then
\begin{align*}
&\lim_{N \to +\infty} \int_{\gamma_{j,N}} (1+|z|)^M e^{\sigma|\Im z|+\tau \Re z} \, dz =0\,,\\
&\lim_{N \to +\infty} \int_{\eta_{j,N}} (1+|z|)^M e^{\sigma|\Im z|+\tau \Re z} \, dz = \int_{-\infty}^{+\infty} (1+|y|)^M e^{\sigma|y|} \; dy < \infty
\end{align*}

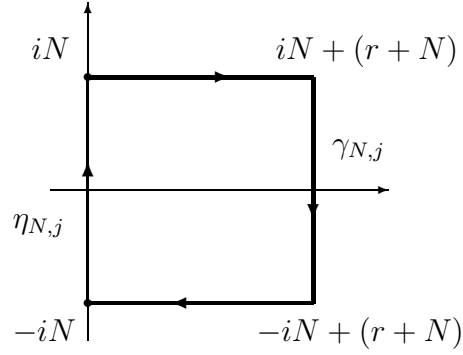
\begin{figure}[h] \label{figure:contourN}
\setlength{\unitlength}{0.5mm}
\begin{picture}(60,80)
\thicklines
\linethickness{0.4mm}
\put(10,70){\vector(1,0){38}}
\put(40,70){\line(1,0){30}}
\put(70,70){\vector(0,-1){38}}
\put(70,40){\line(0,-1){30}}
\put(70,10){\vector(-1,0){38}}
\put(40,10){\line(-1,0){30}}
\linethickness{0.25mm}
\put(10,10){\vector(0,1){38}}
\put(10,40){\line(0,1){30}}
\thinlines
\put(0,40){\vector(1,0){90}}
\put(10,0){\vector(0,1){90}}
\put(-10,0){$-iN$}
\put(-5,75){$iN$}
\put(55,0){$-iN+(r+N)$}
\put(60,75){$iN+(r+N)$}
\put(10,10){\circle*{2}}
\put(10,70){\circle*{2}}
\put(75,50){$\gamma_{N,j}$}
\put(-10,30){$\eta_{N,j}$}
\end{picture}
\caption{$C_{j,N}=\gamma_{j,N} \cup \eta_{j,N}$}
\end{figure}
\end{Lemma}

Write $C_{j,N}=\gamma_{j,N} \cup \eta_{j,N}$ as above.
Then
\begin{align*}
\oint_{C_{1,N}} \cdots \oint_{C_{l,N}}&=\left( \int_{\eta_{1,N}} - \int_{-\gamma_{1,N}}\right)
\cdots \left( \int_{\eta_{1,N}} - \int_{-\gamma_{l,N}}\right)\\
&=\int_{\eta_{1,N}} \cdots \int_{\eta_{l,N}} \pm {\sum_\nu} \int_{\nu_{1,N}} \cdots \int_{\nu_{l,N}}
\end{align*}
where the sum is over all possible combinations $\nu=(\nu_{1,N}, \dots, \nu_{l,N})$,
with $\nu_{j,N} \in \left\{\eta_{j,N}, -\gamma_{j,N}\right\}$, which are different from
$(\eta_{1,N}, \dots, \eta_{l,N})$.
In the following we write $d\lambda$ for  $d\lambda_1 \cdots d\lambda_l$ and assume that the integration is performed in that order.
Then
\begin{multline}
\label{eq:splitting-int}
\oint_{C_{1,N}} \cdots \oint_{C_{l,N}} \frac{a(\l)b(\l)}{c(\l)c(-\l)} \, \varphi_\l(\exp H \cdot o) \;
d\l=
 \int_{\eta_{1,N}} \cdots \int_{\eta_{l,N}}  \frac{a(\l)b(\l)}{c(\l)c(-\l)} \, \varphi_\l(\exp H \cdot o) \;
d\l  \\
\pm {\sum_\nu} \int_{\nu_{1,N}} \cdots \int_{\nu_{l,N}}
\frac{a(\l)b(\l)}{c(\l)c(-\l)} \, \varphi_\l(\exp H \cdot o) \;
d\l\, .
\end{multline}
Here $H \in \overline{\mathfrak a^+}$ is fixed and we suppose that $\|H\| < P/\Omega$. 
Lemma \ref{lemma:integrals} and the estimate (\ref{eq:main-est-shift}) prove that
the first integral on the right-hand side of (\ref{eq:splitting-int}) converges to
$$\int_{i\mathfrak a^*} \frac{a(\l)b(\l)}{c(\l)c(-\l)} \, \varphi_\l(\exp H \cdot o) \; d\l$$
as $N \to +\infty$.
The estimate (\ref{eq:main-est-shift}) also shows that for every $\nu=(\nu_{1,N}, \dots, \nu_{l,N})$
we have
$$
\int_{\eta_{1,N}} \cdots \int_{\eta_{l,N}}  \left|\frac{a(\l)b(\l)}{c(\l)c(-\l)} \, \varphi_\l(\exp H \cdot o) \right|
d\l \leq C_\nu \prod_{j=1}^l \int_{\nu_{j,N}} (1+|\l_j|)^M e^{\sigma|\Im \l_j|+\tau \Re\l_j} \; d\l_j\,,$$
with
$$\sigma=A-\pi <0 \qquad \text{and} \qquad \tau=\|H\|\Omega-P <0\,.$$
According to Lemma \ref{lemma:integrals}, for all $j=1,\dots,l$, we have
$$
\lim_{N\to +\infty} \int_{\nu_{j,N}} (1+|\l_j|)^M e^{\sigma|\Im \l_j|+\tau \Re\l_j} \; d\l_j
=\begin{cases}
0 &\text{if $\nu_{j,N}=\gamma_{j,N}$}\\
\int_{-\infty}^\infty (1+|y|)^M e^{\sigma |y|} \, dy <\infty &\text{if $\nu_{j,N}=\eta_{j,N}$}
\end{cases}
 $$
 Each $\nu=(\nu_{1,N}, \dots, \nu_{l,N})$ contains at least one index $j$ for which
 $\nu_{j,N}=\gamma_{j,N}$.
So
 $$\lim_{N\to +\infty} \int_{\nu_{1,N}} \cdots \int_{\nu_{l,N}} \frac{a(\l)b(\l)}{c(\l)c(-\l)} \, \varphi_\l(\exp H \cdot o) \; d\l=0\,.$$
Thus
$$\lim_{N\to +\infty} \oint_{C_{1,N}} \cdots \oint_{C_{l,N}} \frac{a(\l)b(\l)}{c(\l)c(-\l)} \, \varphi_\l(\exp H \cdot o) \;
d\l=\int_{i\mathfrak a^*} \frac{a(\l)b(\l)}{c(\l)c(-\l)} \, \varphi_\l(\exp H \cdot o) \; d\l\,.$$
By the $W$-invariance of $[c(\l)c(-\l)]^{-1} \, \varphi_\l$, this concludes the proof of Part 2 of Theorem \ref{thm:RamanujanSymm} in the case $\sigma=0$.

To replace the integration along the imaginary axis $i\mathfrak a^*$ by the integration along any
translate $\sigma+i\mathfrak a^*$ with $\sigma\in B(T_\delta)$, it suffices to use (\ref{eq:estWabphilN}) and the following lemma, which is a consequence of Cauchy's theorem.

\begin{Lemma}
Let $D \subset \mathfrak a^*$ be nonempty, compact and connected, and let $F:T_D=D+i\mathfrak a^* \to \C$ be holomorphic. Suppose that, for every compact subset $\omega \subset D$ and every integer $N \geq 0$, we have
$$\sup_{\Re\l\in D} (1+\|\l\|)^N |F(\l)| <\infty\,.$$
Then, for any $\sigma \in D$, the integral $\int_{i\mathfrak a^*} F(\sigma+\l)\, d\l$ exists and is independent of $\sigma$.
\end{Lemma}
\begin{proof}
See e.g. \cite{GV}, Lemma 6.6.2.
\end{proof}

We now prove that the function $f(x)=\int_{i\frak a^*} \wt a(\l) \varphi_\l(x) |c(\l)|^{-2}\, d\l$ extends as a $K_\C$-invariant holomorphic function of $x$ on a neighborhood of $X_G$ in $X_\C$.
Here we have put
$$
\wt a(\l)=\sum_{w \in W}a(w\l)b(w\l)\,.
$$
Recall that $\varphi_\l(x)$ extends holomorphically as a $K_\C$-invariant function on the domain $K_\C \exp(2\Omega_\pi)\cdot o$ with $\Omega_\pi$ given in (\ref{eq:Omegapi}). By $K_\C$-invariance, it therefore suffices to consider the holomorphic extension of $f$ inside $\exp(\Omega_\pi)\cdot o$. On this domain, Opdam's estimates (\ref{eq:OpdamEstimates}) are available.

By estimates (\ref{eq:estimatebcTdelta}) and Remark \ref{rem:equivalentnorms}, there is a constant
$C_\delta>0$ so that
\begin{align*}
|\wt a(\l)| |c(\l)|^{-2}
&\leq C_\delta (1+\|\l\|)^M e^{(A-\pi)\sum_{j=1}^l |\Im\l_j|}\\
&\leq C_\delta (1+\|\l\|)^M e^{(A-\pi)c_2\|\Im\l\|}
\end{align*}
for all $\l \in i\mathfrak a^*$.
For $\varepsilon >0$ we set
$$T_{\varepsilon,\pi} =\{H \in \mathfrak a_\C: \|\Im H\|<\varepsilon\}\,.$$
Suppose we have chosen $\varepsilon$ so that $(\pi -A)c_2 >\varepsilon$ and
$T_{\varepsilon,\pi} \subset \Omega_\pi$.
By (b) of Lemma \ref{lemma:estphil},
there is a constant $C_{\delta,\varepsilon}>0$ such that
$$
|\wt a(\l) \varphi_\l(\exp H \cdot o)| |c(\l)|^{-2} \leq C_{\delta,\varepsilon}
(1+\|\l\|)^M e^{\big((A-\pi)c_2+\varepsilon\big)\|\Im \l\|}\,.
$$
The right-hand side of this inequality is an exponentially decaying function of $\Im \l$,
hence integrable on $Q \times i\mathfrak a^*$ where $Q$ is any compact subset of $T_{\varepsilon,\Omega}$. This allows us to apply to the $f(\exp H \cdot o)$ the theorems of Morera and Fubini, and the claim follows.

Finally, to prove the third part of Ramanujan's Master theorem, we use Lemma \ref{lemma:Schwartztube}, which states that if $a \in \mathcal H(A,P,\delta)$ then $\wt a(\l) \in \mathcal S(\mathfrak a_\varepsilon^*)^W$ for a certain $\varepsilon \in ]0,1[$. This in fact implies that $\mathcal F_G^{-1}\wt a \in
 \mathcal S^p(X_G)^K \subset L^p(X_G)^K \cap L^2(X_G)^K$ with $p=2/(\varepsilon +1) \in ]1,2[$.
Let $f$ be the spherical Fourier series associated with $a$ as in (\ref{eq:fFourier}). Comparison of the inversion formula (\ref{eq:inversionsphericalG}) with (\ref{eq:extensionf-symm}) shows that $\mathcal F_G^{-1}\wt a$ is a smooth $K$-invariant extension of $f$ to all of $G/K$. Formula (\ref{eq:Ramanujansymm}) states then that $\wt a$ is the spherical Fourier transform of $\mathcal F_G^{-1}\wt a$.
The identity holds pointwise for $\l\in T_\varepsilon$ and in  $L^2$ sense for $\l\in i\mathfrak a^*$. The right-hand side of (\ref{eq:Ramanujansymm}) provides a holomorphic extension of the spherical Fourier transform of $\mathcal F_G^{-1}\wt a$ to all of
$T_{\Sigma,m} \cap T_\delta$.
Finally the equality of $L^2$-norms of $f$ and $\wt a$ is an immediate consequence of the Plancherel theorem for $\mathcal F_G$.
This concludes the proof of Theorem \ref{thm:RamanujanSymm}.

\section{The reductive case}
\label{section:RamanujanRed}
\noindent
In this section we extend Ramanujan's Master Theorem to reductive Riemannian symmetric spaces. As in Section \ref{section:prelim}, we consider Riemannian symmetric spaces in duality $X_U=U/K$ and $X_G=G/K$ inside their complexification $X_\C=G_\C/K_\C$. We still assume that $K$ is connected, but we now remove the assumption that
$U$ is semisimple. References for the following structures are  Chapter II in \cite{Take}, Part II \S 1 in \cite{Varadarajan},
and Sections 1 and 2 in \cite{BOP}.

Let $\mathfrak z$ be the center of the Lie algebra $\mathfrak u$ of $U$. Then $\mathfrak u=\mathfrak z \oplus \mathfrak u'$ where $\mathfrak u'=[\mathfrak u,\mathfrak u]$ is semisimple. As in Section
\ref{section:prelim}, let $\tau$ be the involution associated with $X_U$, and let $\mathfrak u=\mathfrak k \oplus i\mathfrak p$ be the corresponding decomposition of $\mathfrak u$. Notice that $\tau$ preserves $\mathfrak z$ and $\mathfrak u'$. 
We shall assume that $\mathfrak k \cap \mathfrak z=\{0\}$, i.e. that the symmetric pair $(\mathfrak u,\mathfrak k)$ is effective. Hence
$\mathfrak u'=\mathfrak k \oplus i\mathfrak p'$ with $\mathfrak p'=\mathfrak p \cap (i\mathfrak u')$  and $i\mathfrak p=\mathfrak z \oplus (i\mathfrak p')$.
The Lie algebra $\mathfrak g= \mathfrak k\oplus \mathfrak p$ of $G$ is reductive, and $\mathfrak g=i\mathfrak z \oplus \mathfrak g'$ where $i\mathfrak z$ is the center of $\mathfrak g$ and $\mathfrak g'=[\mathfrak g,\mathfrak g]=\mathfrak k\oplus \mathfrak p'$ is semisimple.

Set $\Gamma_0=\{X \in \mathfrak z:\exp X=e\}$ where $e$ is the identity of $U$. Then $\Gamma_0$ is a full rank lattice in $\mathfrak z$ and $T=\mathfrak z/\Gamma_0=\exp \mathfrak z$ is isomorphic to the identity component of the center of $U$.   Let $U'$ be the analytic subgroup of $U$ with Lie algebra $\mathfrak u'$. Then $U'$ is a compact connected semisimple
Lie group with finite center containing $K$. Moreover $U=TU' \cong T \times_F U'$  where $F=T \cap U'$ is a finite central subgroup of $U$.
We shall assume for simplicity that $F$ is trivial. Hence $U \cong T \times U'$.

The involutive automorphism $\tau$ leaves $K$ invariant and $X_{U'}=U'/K$ is a semisimple Riemannian symmetric space of the compact type. Moreover
\begin{equation} \label{eq:UKred}
U/K \cong T \times U'/K \cong  \exp(\mathfrak z) \times U'/K\,.
\end{equation}
Let $G'$ be the analytic subgroup of $G$ with Lie algebra $\mathfrak g'$. Then $G'$ is a noncompact connected semisimple Lie group with finite center, and $K$ is maximal compact in $G'$. The subgroup $V=\exp(i\mathfrak z)$ is the split component of $G$. We have $G=G'V$ with $G' \cap V=\{e\}$. Hence
\begin{equation} \label{eq:GKred}
G/K \cong V  \times G'/K= \exp(i\mathfrak z) \times G'/K\,.
\end{equation}
Furthermore, $G_\C=(G'V)_\C \cong V_\C \times G'_\C$, where $V_\C=\exp(\mathfrak z \oplus
i\mathfrak z)$. So
\begin{equation} \label{eq:GcKcred}
G_ \C/K_\C \cong V_\C  \times G'_\C/K_\C\,.
\end{equation}

Let  $\mathfrak a  \subset \mathfrak p$ be a maximal abelian subspace. Then $\mathfrak a=i\mathfrak z \oplus
\mathfrak a'$ where $\mathfrak a' \subset \mathfrak p'$ is maximal abelian.
We fix an inner product $\inner{\cdot}{\cdot}$ on $\mathfrak a$ by setting it on $\mathfrak a'$ equal to the one associated with the Killing form, equal on $i\mathfrak z$ to a fixed inner product, and by declaring that $\mathfrak a'$ and  $i\mathfrak z$ are orthogonal in $\mathfrak a$.  Extend then $\inner{\cdot}{\cdot}$ by duality on $\mathfrak a^*$ and by $\C$-bilinearity on $\mathfrak a_\C$ and $\mathfrak a^*_\C$. We denote by $\|\cdot\|$ the norms on $\mathfrak a$ and $\mathfrak a^*$ associated with the inner product $\inner{\cdot}{\cdot}$.
If $\l \in \mathfrak a_\C^*$ and $a=\exp H \in \exp(\mathfrak a_\C) \subset G_\C$, then we write
$a^\l=e^{\l(H)}$, provided this is well defined.

Let $\Sigma$ be the set of restricted roots of $(\mathfrak g',\mathfrak a')$. Denote by $\Sigma^+$ a choice of positive roots and by $(\mathfrak a')^+$ the corresponding positive Weyl chamber.
Finally set $\mathfrak a^+=i\mathfrak z \oplus (\mathfrak a')^+$. The Weyl group $W$ of $(\mathfrak g,\mathfrak a)$ is the finite group generated by the reflections relative to $\Sigma$. It acts trivially on  $i\mathfrak z$. A fundamental domain for the action of $W$ on $\mathfrak a$ is $\overline{\mathfrak a^+}=i\mathfrak z \oplus \overline{(\mathfrak a')^+}$.

The definitions of spherical functions on semisimple Riemannian symmetric spaces via (\ref{eq:varphil}) and (\ref{eq:psimu})
extend to the reductive case.  Notice also that, by Remark \ref{rem:notsimplyconnected}, we can
suppose that $U'/K$ is simply connected.

As proven in \cite{BOP}, Section 2, the set parametrizing the
$K$-spherical representations of $U$ is in this case $\Lambda^+_K(U)=i\Gamma_0^* \oplus \Lambda^+_K(U')$ where
\begin{equation}
\label{eq:Gammastar}
i\Gamma_0^*=\{\mu^0 \in i\mathfrak z^*: \text{$\mu^0(H)\in 2\pi  i\Z$ for all $H \in \Gamma_0$}\}
\end{equation}
and $\Lambda^+_K(U')= \{\mu\in(\mathfrak a')^*: \text{$\mu_\a\in\Z^+$ for all $\a\in\Sigma^+$}\}$ is as in
(\ref{eq:LambdaK}). Here the direct sum symbol means that every element $\mu \in \Lambda^+_K(U)$ admits a unique decomposition as a sum of an element of $\mu^0 \in i\Gamma_0^*$ and an element of $\mu'\in \Lambda^+_K(U')$.
If $\mu=\mu^0 +\mu'$, then $d(\mu)=d(\mu')$. Moreover, the spherical function of spectral parameter $\mu=\mu^0 +\mu'$ is
\begin{equation}
\label{eq:psimured}
\psi_\mu(tu')=t^{\mu^0} \psi_{\mu'}(u')\,, \qquad t \in T=\exp(\mathfrak z), \, u' \in U'\,,
\end{equation}
where $\psi_{\mu'}(u')$ is the spherical function of spectral parameter $\mu'$ on the semisimple Riemannian symmetric space of the compact type $U'/K$. A similar property holds for the spherical functions on $G/K$: if $\l=\l^0+\l' \in
\mathfrak a_\C^* =\mathfrak z_\C^* \oplus (\mathfrak a')^*_\C$, then the spherical function
$\varphi_\l$ on $G/K$ is given by
\begin{equation}
\label{eq:varphilred}
\varphi_\l(xg')=x^{\l^0} \varphi_{\l'}(g')\,, \qquad x \in V=\exp(i\mathfrak z), \, g' \in G'\,,
\end{equation}
where $\varphi_{\l'}(g')$ is the spherical function of spectral parameter $\l'$ on the semisimple Riemannian symmetric space of the noncompact type $G'/K$. In particular, as in the semisimple case, the spherical functions on $U/K$ extend holomorphically to $G_\C$ and
\begin{equation}
\psi_{\mu}|_G=\varphi_{\mu+\rho}
\end{equation}
where $\rho=1/2 \sum_{\a \in \Sigma^+} m_\a \a \in (\mathfrak a')^*$.
Since $\inner{\l}{\a}=0$ for $\l\in\mathfrak z^*_\C$ and $\a\in\Sigma$, we can extend the definition of the $c$-function to
$\mathfrak a^*_\C$ by the same formula (\ref{eq:c}) as in the semisimple case. We then have $c(\l)=c(\l')$ if
$\l=\l^0+\l' \in\mathfrak z^*_\C \oplus (\mathfrak a')^*_\C$.

Formulas (\ref{eq:psimured}) and (\ref{eq:varphilred}) reduce the spherical harmonic analysis on the pair of
reductive symmetric spaces $U/K$, $G/K$ to the harmonic analysis on the abelian spaces $T$, $V$  together
with the spherical harmonic analysis on the semisimple symmetric spaces $U'/K$, $G'/K$.

Let $v=\dim V$ be the dimension of the split component of $G$. Since $\Gamma_0$ is a full rank lattice in $\mathfrak z$ we can choose linearly independent vectors $e_1, \dots, e_v \in \mathfrak z$ so that $\Gamma_0=\sum_{k=1}^v \Z e_k$.
Define $\varepsilon_1, \dots, \varepsilon_v \in i\Gamma_0^*$ by $\varepsilon_k(e_h)=2\pi i \delta_{k,h}$. Then  $\{\varepsilon_1, \dots, \varepsilon_v\}$ is a basis of $(i\mathfrak z)^*$.
We fix $\{\varepsilon_1, \dots, \varepsilon_v, \omega_1, \dots, \omega_l\}$ as a basis of
$\mathfrak a^*=(i\mathfrak z)^* \oplus (\mathfrak a')^*$.
The corresponding decomposition of $\l=\l^0+\l' \in\mathfrak a_\C^*=\mathfrak z_\C^* \oplus (\mathfrak a')_\C^*$ will be written
either as 
\begin{equation} \label{eq:lambda-coords-reductive}
\l=\sum_{j=1}^v \l_j \varepsilon_j+\sum_{j=v+1}^{v+l} \l_j \omega_j \qquad\text{or as} \qquad
\l=\sum_{k=1}^v \l^0_k \varepsilon_k+\sum_{j=1}^{l} \l'_j \omega_j\,. 
\end{equation}
Define
\begin{align}
\Gamma_0^+&=\sum_{k=1}^v \Z^+ e_k\,,\\
\label{eq:Gammastarplus}
(i\Gamma_0^*)^+&=\{\mu^0 \in i\mathfrak z^*: \text{$\mu^0(H)\in 2\pi  i\Z^+$ for all $H \in \Gamma_0^+$}\}
\end{align}
and
\begin{equation}
\label{eq:Lambdaplusplus}
\Lambda^{++}=(i\Gamma_0^*)^+ \oplus \Lambda^+_K(U')\,.
\end{equation}
If $\mu=\mu^0+\mu' \in \Lambda^{++}$, then $\mu^0_k \in \Z^+$ for all $k=1,\dots, v$ and $\mu'_j \in \Z^+$ for all $j=1,\dots, l$.
We set
$$
|\mu|=\mu_1^0 +\dots +\mu_v^0+\mu'_1+\dots+\mu'_l\,.
$$
Moreover, if $x=\exp X$ with $X=i\sum_k X_k e_k \in i \mathfrak z$ and $\l^0=\sum_k \l^0_k \varepsilon_k \in
\mathfrak z^*_\C$, then
$$x^{\l ^0}=\prod_{k=1}^v e^{\l^0_k X_k}=\prod_{k=1}^v x_k^{\l^0_k} \qquad \text{with} \qquad x_k=e^{ X_k} \in ]0,+\infty[\,.$$

Let $A, P, \delta$ be constants so that $A<\pi$, $P>0$ and $0<\delta \leq 1$.
Set
\begin{equation}
\label{eq:Hdeltared}
\mathcal H(\delta)=\left\{\l=\l^0+\l' \in \mathfrak a^*_\C: \; \text{$\Re\l'_\b > -\delta \widetilde\rho_\b$ for $\b \in \Sigma^+_*$ and
$\Re\l^0_k>-\delta$ for $k=1,\dots, v$}\right\}\,.
\end{equation}
The Hardy class $\mathcal H(A,P,\delta)$ is the space of all holomorphic functions $a:\mathcal H(\delta) \to \C$ satisfying the  growth condition:
there exists a constant $C>0$ so that
\begin{equation}
 \label{eq:Hardyclassred}
|a(\l)|\leq C \prod_{j=1}^{v+l} e^{-P(\Re\l_j)+A|\Im \l_j|}
\end{equation}
for all $\l \in \mathcal H(\delta)$.
Furthermore, set
\begin{equation} \label{eq:bred}
 b(\l)=b^0(\l^0)b'(\l')\,, \qquad \l=\l^0+\l' \in \mathfrak a_\C^*\,,
\end{equation}
where $b'(\l')$ is the function $b$ of the semisimple Riemannian symmetric spaces $U'/K$, $G'/K$, as in (\ref{eq:b}), and
\begin{equation}
 \label{eq:b0}
b^0(\l^0)=\left(\frac{i}{2}\right)^v \prod_{k=1}^v \frac{1}{\sin(\pi \l^0_k)}\,.
\end{equation}
Observe that, since $W$ acts trivially on $\mathfrak z^*_\C$, we have
\begin{equation}
\sum_{w \in W} a(w\l)b(w\l)=b^0(\l^0) \sum_{w \in W} a(w\l)b'(w\l')\,, \qquad \l=\l^0+\l'\,.
\end{equation}
Let $T'_\delta$ be the tube domain in $(\mathfrak a')^*_\C$ defined by  (\ref{eq:Trho}), and let
$$
T^0_\delta=\{\l^0 \in \mathfrak z^*_\C: \text{$0<\Re\l^0_k <\delta$ for all $k=1,\dots,v$}\}.
$$
Consider the tube domain in $\mathfrak a^*_\C$ given by
\begin{equation}
\label{eq:Trhored}
T_\delta=T^0_\delta \oplus T'_\delta\,.
\end{equation}
Its base in $\mathfrak a^*$ is $B(T_\delta)=\{\l^0 \in (i\mathfrak z)^*: \text{$0<\l^0_k <\delta$ for all $k=1,\dots,v$}\} \oplus B(T_\delta)$.
Finally, let $T'_{\Sigma,m}$ be as in (\ref{eq:TSigmam}) and let $T_{\Sigma,m}=\mathfrak z^*_\C \oplus T'_{\Sigma,m}$.

The following theorem combines the semisimple and the multivariable abelian versions of Ramanujan's Master Theorem.
Because of formulas (\ref{eq:lambda-coords-reductive}) to (\ref{eq:Trhored}),  
its proof reduces to a straightforward combination of the arguments used in these two cases, and is omitted.

\begin{Thm}[Ramanujan's Master Theorem for reductive Riemannian symmetric spaces]
\label{thm:RamanujanSymmRed}

Keep the above assumptions and notation, and let $a \in \mathcal H(A,P,\delta)$. Then the following properties hold:
\begin{enumerate}
\item
The spherical Fourier series
\begin{equation} \label{eq:fFourier-red}
f(x)=\sum_{\mu\in \Lambda^{++}} (-1)^{|\mu|} d(\mu) a(\mu+\rho) \psi_\mu(x)
\end{equation}
converges normally on  compact subsets of $D_{P/\Omega}=U \exp\big((P/\Omega)B\big) \cdot o_\C$
where $B=\{H \in \mathfrak a: \|H\|<1\}$ is the open unit ball in $\mathfrak a$.
Its sum is a $K$-invariant holomorphic function on the neighborhood $D_{P/\Omega}$ of $X_U$ in $X_\C$.
\item Let $T_\delta$ be the tube domain in (\ref{eq:Trhored}) and let  $\sigma\in B(T_\delta)$. Then for $x=\exp H \in A$ with $\|H\|<P/\Omega$, we have
\begin{equation}\label{eq:extensionf-symm-red}
f(x)=\frac{1}{|W|} \int_{\sigma+i\mathfrak a^*} \left( \sum_{w \in W} a(w\l)b(w\l)\right) \varphi_\l(x) \; \frac{d\l}{c(\l)c(-\l)} \,.
\end{equation}
\item
The formula
\begin{equation} \label{eq:Ramanujansymm-red}
\int_{X_G} f(x)\varphi_{-\l}(x)\; dx=  \sum_{w \in W} a(w\l)b(w\l)
\end{equation}
holds for the extension of $f$ to $X_G$ and for all
$\l \in T_\delta \cap T_{\Sigma,m}=T^0_\delta \oplus T'_{\Sigma,m}$.
\end{enumerate}
\end{Thm}

\section{Further remarks and open problems}
\label{section:further}
\noindent
In \cite{Bertram}, Bertram presents an interesting alternative approach to Ramanujan's Master
Theorem for the rank-one case, by means of the kernel function $k(u,w)=\big(z(u)+z(w)\big)^{-2\rho}$. In this formula,
$z$ is the holomorphic continuation to $X_\C$ of the function defined on $X_U$ by $z(x)=\cos(d(x,K))$, where $d(x,K)$ is
the distance of the point $x$ from the base point in $X_U$.
Moreover, $2\rho=m_{\b/2}/2+m_\b$ for the unique root $\b\in\Sigma_*^+$.
The function $z$ is natural, as every $K$-invariant function on $X_U$ factors through $z$.
By considering $k(x,y)=k_x(y)=k_y(x)$, Bertram proves analogues to Mehler's and Neumann's formulas, stating that there
exists meromorphic functions $d(\l)=d(-\l)$ and $e(\l)$ (explicitly given as ratios of products of gamma functions) so that
$$(\mathcal F_Gk_y)(\l)=d(\l)\varphi_\l(y) \quad \text{for}\quad  |\Re\l_\b|<\rho$$
and
$$(\mathcal F_U k_x)(\mu+\rho)=(-1)^\mu e(\mu+\rho) \Phi_{-\mu-\rho}(x) \quad \text{for}\quad \mu\in\Z^+\,.$$
Here $\Phi_\l$ denotes Harish-Chandra's series. By means of Harish-Chandra's relation
$\varphi_\l=c(\l)\Phi_\l+c(-\l)\Phi_{-\l}$, he finds $b(\l)=\frac{d(\l)}{e(\l)}\, c(-\l)$. The explicit formulas
for $d, e$ and $c$  allow him to recover the formula for $b$. It would be interesting to find a similar approach in
higher rank as well. Unfortunately, we do not know higher dimensional analogues neither of the kernel function nor
of the formulas by Mehler or Neumann, and the knowledge of $b$ provides information only on the ratio $\frac{d(\l)}{e(\l)}$.

Another open question is on the nature of the functions on $X_G$, $X_U$ and $X_\C$ satisfying the assumptions stated in Ramanujan's theorem. For instance, what are the functions on $X_U$ having Fourier coefficients coming from elements of $\mathcal H(A,P,\delta)$? This kind of questions are related to theorems of Paley-Wiener type. Notice that one can use Laplace transforms to obtain elements of
$\mathcal H(A,P,\delta)$. Indeed, identify $\mathfrak a_\C^*$ with $\C^l$ by means of the basis $\Pi_*$,
as in section \ref{subsection:coordinates}.
Let $h:\mathfrak a_\C^*\equiv \C^l \to \C$ be integrable and having support inside a domain of the form $[P,R]^l$ with
$P<R<+\infty$, and set
$$\mathcal Lh(\l)=\int_0^{+\infty}\!\!\!\dots\!\int_0^{+\infty} h(x) e^{-\sum_{j=1}^l \l_j x_j} \; dx\,.$$
Then $\mathcal Lh(\l)$ is an entire function on $\mathfrak a_\C^*$. Moreover, for every $d>0$ there is a constant $C_d>0$ so
that $|(\mathcal Lh)(\l)|\leq C_d e^{-P\sum_{j=1}^l \Re \l_j}$ for
all $\l\in \mathfrak a_\C^*$ with $\Re\l_j>-d$ for all $j$. Thus $\mathcal Lh \in \mathcal H(A,P,\delta)$ for all $0\leq A <\pi$.

Finally, the spherical Fourier transform on Riemannian symmetric spaces of the noncompact and of the compact type
has been extended by the works of Heckman, Opdam and Cherednik to the setting of hypergeometric functions
associated with root systems.
See \cite{OpdamActa}, \cite{OpdamJapon} and references therein.
A natural question is therefore a generalization of Ramanujan's Master Theorem in this setting.
The necessary $L^p$-harmonic analysis on root systems needed for instance to generalize the proof of the final part of Theorem \ref{thm:RamanujanSymm} has been recently developed in \cite{NPP}. There are nevertheless technical difficulties,
for instance the fact that at some points we use the classification of the root multipliticies,
listed at the beginning of section \ref{section:PlancherelDensities}.
These kinds of arguments have to be replaced. 
A root system situation where all problems should be easily overpassed corresponds to the even multiplicity case.
Indeed, the shift operators techniques developed by Opdam allow in this case a reduction to a Euclidean
$W$-invariant situation where a multivariable analogue of the classical Ramanujan's theorem can be easily proven.
See \cite{HS94}, Part 1, and  \cite{OPeven,OPeven1}.
We will come back to these issues in future work.

\end{document}